\documentclass[11pt]{article}

\usepackage{latexsym,amsfonts,amssymb,amsthm}
\usepackage[centertags]{amsmath}

\allowdisplaybreaks

\theoremstyle{plain}
\newtheorem{mtheorem}{Main Theorem}

\swapnumbers
\theoremstyle{plain}
\newtheorem{theorem}{Theorem}[section]
\newtheorem{proposition}[theorem]{Proposition}
\newtheorem{corollary}[theorem]{Corollary}
\newtheorem{lemma}[theorem]{Lemma}
\newtheorem*{theorem*}{Theorem}
\newtheorem{remark}[theorem]{Remark}

\begin{document}

%
%

\newcounter{ListingNr}

\newenvironment{listing}%
{\begin{list}   { {\rm (\roman{ListingNr})} }
                { \topsep2mm
                  \partopsep0ex
                  \listparindent0pt
                  \itemsep0.5ex
                  \parskip0pt
                  \leftmargin5ex
                  \labelwidth5.5ex
                  \labelsep0ex
                  \renewcommand{\makelabel}[1]{\textit{##1}\hfill}
                  \usecounter{ListingNr}} }%
{\end{list}}

\newcommand{\mlt}[3]{[#1,#2]\raisebox{-3pt}{${}_{#3}$}}
\newcommand{\pairing}[2]{\langle\,#1,#2\,\rangle}

\newcommand{\dank}{\unitlength1pt%
  \begin{picture}(8,8)%
    \linethickness{0.4pt}
    \put(1,0){$\downarrow$}%
    \put(1.1,0.7){\line(1,0){5}}
  \end{picture}}

\newcommand{\upnk}{\unitlength1pt%
  \begin{picture}(8,8)%
    \linethickness{0.4pt}
    \put(1,0){$\uparrow$}%
    \put(1.15,4.8){\line(1,0){4.8}}
  \end{picture}}

\newcommand{\inn}{\unitlength1pt%
  \begin{picture}(7,8)%
    \put(1.8,1.2){{\tiny $*$}}%
  \end{picture}}

\newcommand{\descn}{\unitlength1pt%
  \raisebox{-6pt}{%
    \begin{picture}(22,9)%
      \put(4,6){$\mathop{\smile}$}%
      \put(11,0){$\scriptstyle D$}
    \end{picture}}
  }
\newcommand{\descr}{\unitlength1pt%
  \raisebox{-6pt}{%
    \begin{picture}(21,8)%
      \put(4,6){$\mathop{\sim}$}%
      \put(10,0){$\scriptstyle D$}
    \end{picture}}
  }

\newcommand{\peakn}{\unitlength1pt%
  \raisebox{-6pt}{%
    \begin{picture}(22,9)%
      \put(4,6){$\mathop{\smile}$}%
      \put(11,0){$\scriptstyle P$}
    \end{picture}}
  }
\newcommand{\peaknumn}{\unitlength1pt%
  \raisebox{-6pt}{%
    \begin{picture}(24,9)%
      \put(4,6){$\mathop{\smile}$}%
      \put(11,0){$\scriptstyle |P|$}
    \end{picture}}
  }
\newcommand{\speakn}{\unitlength1pt%
  \raisebox{-6pt}{%
    \begin{picture}(22,9)%
      \put(4,6){$\mathop{\smile}$}%
      \put(11,2){$\scriptstyle s^*$}
    \end{picture}}
  }
\newcommand{\peakr}{\unitlength1pt%
  \raisebox{-6pt}{%
    \begin{picture}(21,8)%
      \put(4,6){$\mathop{\sim}$}%
      \put(10,0){$\scriptstyle P$}
    \end{picture}}
  }
\newcommand{\peaknumr}{\unitlength1pt%
  \raisebox{-6pt}{%
    \begin{picture}(23,8)%
      \put(4,6){$\mathop{\sim}$}%
      \put(10,0){$\scriptstyle |P|$}
    \end{picture}}
  }

\newcommand{\breakitdown}{\begin{center} $*$ $*$ $*$ \end{center}}

\newcommand{\assozerl}{%
                       \unitlength1pt
                       \raisebox{-2pt}{
                         \begin{picture}(9,10)
                           \put(1,0){\line(0,1){10}}
                           \put(1,2){\mbox{$\ass$}}
                         \end{picture}
                       }
                      }

\newcommand{\sctriangle}{{\scriptstyle\triangle}}
\newcommand{\scsctriangle}{{\scriptscriptstyle\triangle}}
\newcommand{\Cl}{\mathsf{Cl}}
 
\newcommand{\C}{\mathbb{C}}
\newcommand{\N}{\mathbb{N}}

\newcommand{\DD}{\mathcal{D}}
\newcommand{\LL}{\mathcal{L}}
\newcommand{\EE}{\mathcal{E}}
\newcommand{\ZZ}{\mathcal{Z}}
\newcommand{\KK}{\mathcal{K}}
\newcommand{\PP}{\mathcal{P}}
\renewcommand{\AA}{\mathcal{A}}
\newcommand{\PPP}{\mathfrak{P}}
\newcommand{\HH}{\tilde{\mathcal{H}}}
\newcommand{\CC}{\mathcal{C}}
\newcommand{\maj}{\mathrm{maj}\,}
\newcommand{\minimaj}{\mathrm{maj}\,}

\newcommand{\Dbar}{\overline{D}}
\renewcommand{\P}{\mathsf{P}}
\newcommand{\T}{\mathsf{T}}
\newcommand{\LSC}{\mathsf{LSC}}
\newcommand{\Qsym}{\mathsf{Qsym}}
\newcommand{\Part}{\mathsf{Part}}
\newcommand{\Prim}{\mathsf{Prim}}
\newcommand{\cha}{\mathsf{char}}
\newcommand{\ch}{\mathsf{ch}}
\newcommand{\length}{\ell}
\newcommand{\Des}{\mathsf{Des}}
\newcommand{\SP}{\mathsf{SP}}
\newcommand{\Peak}{\mathsf{Peak}}
\newcommand{\peak}{\mathsf{peak}}
\newcommand{\sign}{\mathsf{sign}}
\newcommand{\ord}{\mathsf{ord}}
\newcommand{\Rad}{\mathsf{Rad}}
\renewcommand{\dim}{\mathsf{dim}}
\renewcommand{\ker}{\mathsf{ker}}
\newcommand{\R}{\tilde{\Xi}}

\newcommand{\eps}{\varepsilon}
\newcommand{\da}{\downarrow}
\newcommand{\ass}{\approx}
\newcommand{\conv}{\star}

\newcommand{\summe}{\mathsf{sum}\,}
\newcommand{\ol}[1]{\overline{#1}}

\newcommand{\Kasten}[1]{
                       \begin{picture}(1,1)
                        \put(0,0){\makebox(1,1){#1}}
                        \put(0,0){\line( 1, 0){1}}
                        \put(0,0){\line( 0, 1){1}}
                        \put(1,1){\line(-1, 0){1}}
                        \put(1,1){\line( 0,-1){1}}
                       \end{picture}}

\newcommand{\tensor}{\otimes}
\newcommand{\emptyword}{\varnothing}

\def\haken#1{\underline{#1}%
             {\raise -0.35ex%
              \hbox{\vphantom{$#1$}\vrule height 0.8ex}}}

\newcommand{\lra}{\longrightarrow}
\newcommand{\lmt}{\longmapsto}
\newcommand{\p}{\vdash}
\newcommand{\zerl}{\models}
\newcommand{\zerlodd}{\models_{\mbox{\tiny odd}}}
\newcommand{\podd}{\vdash_{\mbox{\tiny odd}}}

\newcommand{\Q}{\mathbb{Q}} 
\newcommand{\id}{\mathsf{id}}

\newcommand{\setofall}[2]{\mbox{$\{\,#1\,|\,#2\,\}$}}
\newcommand{\Setofall}[2]{\mbox{$\Big\{\,#1\,\Big|\,#2\,\Big\}$}}
\newcommand{\erz}[1]{\langle\,#1\,\rangle}
\newcommand{\Erz}[1]{\Big\langle\,#1\,\Big\rangle}
\newcommand{\tiptop}[4]{\left\{\begin{array}{ll}
      #1 & #2\\[1mm]
      #3 & #4
    \end{array}\right.}        
\newcommand{\kref}[1]{(\ref{#1})}

\newcommand{\wldots}{.\,...\,.}

%
%

%
%
%
%
%
%
%
%
%
%
%
%
%
%
%
%
%
%
%
%
%
%
%
%

\title{The peak algebra of the symmetric group revisited}

\author{%
  Manfred Schocker\thanks{supported by the Research Chairs of Canada.}\\[2mm]
   \begin{minipage}{5.5cm}
     \begin{center}
          \begin{small}
            LaCIM\\
            Universit{\'e} du Qu{\'e}bec \`a Montr{\'e}al\\
            CP 8888, succ. Centre-Ville\\
            Montr{\'e}al (Qu{\'e}bec), H3C 3P8\\
            Canada\\
            E-mail: \texttt{mschock@math.uqam.ca}
          \end{small}
     \end{center}
   \end{minipage}
  }

\maketitle

\begin{center}
\end{center}

\newpage

%
%
%
%
%
%
%
%
%
%
%
%
%
%
%
%
%
%
%
%
%
%
%
%

\begin{abstract} \noindent
  The linear span $\PP_n$ of the sums of all permutations in the
  symmetric group $S_n$ with a given set of peaks is a sub-algebra of
  the symmetric group algebra, due to Nyman. This peak algebra is a
  left ideal of the descent algebra $\DD_n$; and the direct sum $\PP$
  of all $\PP_n$ is a Hopf sub-algebra of the direct sum $\DD$ of all
  $\DD_n$, dual to the Stembridge algebra of peak functions.  In our
  self-contained approach, peak counterparts of several results on the
  descent algebra are established, including a simple combinatorial
  characterization of the algebra $\PP_n$; an algebraic
  characterization of $\PP_n$ based on the action on the
  Poincar\'e-Birkhoff-Witt basis of the free associative algebra; the
  display of peak variants of the classical Lie idempotents; an
  Eulerian-type sub-algebra of $\PP_n$; a description of the Jacobson
  radical of $\PP_n$ and its nil-potency index, of the principal
  indecomposable and irreducible $\PP_n$-modules, and of the Cartan
  matrix of $\PP_n$.  Furthermore, it is shown that the primitive Lie
  algebra of $\PP$ is free, and that $\PP$ is its enveloping algebra.
\end{abstract}

\newpage\

%
%
%
%
%
%
%
%
%
%
%
%
%
%
%
%
%
%
%
%
%
%
%
%
%

\section{Introduction} \label{intro}

For any positive integer $n$ and any permutation $\pi$ in the
symmetric group $S_n$,
$$
\Des(\pi) := \setofall{i}{1\le i\le n-1,\,i\pi>(i+1)\pi}
$$
is the \emph{descent set} of $\pi$.  In a remarkable paper of 1976,
Solomon discovered (in the wider context of finite Coxeter groups)
that the linear span $\DD_n$ of the elements
$$
\Delta^D := \sum_{\Des(\pi)=D} \pi, \qquad D\subseteq
\{1,\ldots,n-1\}
$$
is a sub-algebra of the group ring $KS_n$ (\cite{solomon76}) of
dimension $2^{n-1}$. The ground ring $K$ will throughout be assumed to
be a field of characteristic~$0$.  Several independent proofs of
Solomon's result have been given, adding up to a better understanding
of the \emph{descent phenomenon} (see, for instance,
\cite{blelau93a,garsia-reutenauer89,gessel84,willigenburg98}).

Malvenuto and Reutenauer introduced the structure of a graded
Hopf algebra on the direct sum
$
KS
=
\bigoplus_{n\ge 0} KS_n
$,
based on the convolution, or outer product
(\cite{malvenuto-reutenauer95}).
The Solomon descent algebra 
$$
\DD
=
\bigoplus_{n\ge 0} \DD_n
$$
is a Hopf sub-algebra of $KS$, dual to the algebra $\Qsym$ of
quasi-symmetric functions (\cite{gessel84,malvenuto-reutenauer95}) and
isomorphic to the algebra of noncommutative symmetric functions
(\cite{parisI}).  The algebra $\DD$ with its triple algebraic
structure has been subject to rapidly increasing interest for the past
fifteen years.  As a consequence, a number of surprising results on
$\DD$ have been obtained, connecting the descent algebra to many other
fields of algebraic combinatorics such as the theory of the free Lie
algebra, the representation theory of the symmetric group, and the
theory of $0$-Hecke algebras and related structures.  Some of the
major contributions are
\cite{atkinson92,blelau96,garsia-reutenauer89,gessel-reutenauer93,%
  joellenbeck99,joellereut01,parisI,reutenauer93}.  Another interesting
branch of the investigations dealt with extensions of the descent
algebra (\cite{loday-ronco98,patrasreut01,schocker-lia}).

\breakitdown

For any permutation $\pi\in S_n$,
$$
\Peak(\pi) := \setofall{i}{2\le i\le n-1,\,(i-1)\pi<i\pi>(i+1)\pi}
$$
is the \emph{peak set} of $\pi$. Our starting point is the
following refinement of \cite{nyman02} (see also \cite{abn02}).

\begin{mtheorem} \label{left-ideal}
  The linear span $\PP_n$ of the
  elements
  $$
  \Pi^P := \sum_{\Peak(\pi)=P} \pi, \qquad P\subseteq
  \{2,\ldots,n-1\}
  $$
  is a left ideal of $\DD_n$ of dimension $f_n$, for all $n\ge 0$,
  where $f_n$ denotes the $n$-th Fibonacci number, defined by
  $f_0=f_1=f_2=1$ and $ f_n=f_{n-1}+f_{n-2} $ for all $n>2$.
\end{mtheorem}
(To save trouble, let it be said that $\PP_n$ is an algebra
\emph{without} identity.)

In addressing the properties of the peak algebra $\PP_n$, numerous
combinatorial and algebraic concepts related to the descent algebra
$\DD_n$ deserve to be transferred to the peak setting.  To emphasize
this, is the main goal of the present work.  Almost every aspect of
the theory of the Solomon descent algebra has a peak counterpart
indeed.  We should like to point out that a new and illuminating
approach has been introduced in \cite{abn02} only recently, relating
$\PP_n$ to the descent algebras of type $B$ and $D$.

As one part of our investigations, three independent proofs of Main
Theorem~\ref{left-ideal} are given, providing a better understanding
of the \emph{peak phenomenon}. The first one is based on the fact that
the direct sum
$$
\PP=\bigoplus_{n\ge 0}\PP_n
$$
is a Hopf sub-algebra of $\DD$ (Section~\ref{peak}); the second one
is purely combinatorial (Section~\ref{comb-approach}); while the third
one is immediate from an algebraic characterization of $\PP_n$ arising
from the action on the Poincar\'e-Birkhoff-Witt basis of the free
associative algebra (Section~\ref{alg-approach}).  

The Hopf algebra $\PP$ is the graded dual of the Stembridge algebra of
peak functions --- a sub-algebra of $\Qsym$, which arose naturally
from the context of enriched $P$-partitions (\cite{stembridge97}).
Stembridge's peak algebra has already been studied intensively
(\cite{bmsw02,bhw02,stembridge97}).  Some of those results are
recovered. The approach presented here, however, is self-contained and
independent of the theory of quasi-symmetric functions; for the
convenience of the reader more familiar with that setting, enough
information is provided to transfer the results to $\Qsym$, by
duality.

\breakitdown

This paper is subdivided into several sections, each of which contains
one main result. Their labeling thus coincides with the labeling of
the sections.  For a first impression, taking note just of these main
results is sufficient. A more detailed overview follows below.  We
mention that the fruits of
\cite{blelau96,blelau02,garsia-reutenauer89,malvenuto-reutenauer95}
will be of particular importance.

In Section~\ref{descent}, some basic results on $\DD$ and its duality
to $\Qsym$ are recalled. This allows to confirm that $\PP$ --- as a
co-algebra --- is the graded dual of Stembridge's peak algebra.
Independent from that, the Hopf algebra structure on $\PP$ is
established in Section~\ref{peak}, which yields, in particular, a
first proof of Main Theorem~\ref{left-ideal}, and recovers results due
to Stembridge \cite{stembridge97} and N. Bergeron et al.
\cite{bmsw02}, by duality.  Furthermore, $\PP$ is a polynomial ring in
the (non-commuting) variables $\Gamma^n$, $n$ odd, with respect to the
outer product, where $\Gamma^n$ denotes the sum of all
\emph{alternating} permutations in $S_n$.

A simple combinatorial characterization of $\PP_n$ follows in
Section~\ref{comb-approach}, which says that $\varphi\in\DD_n$ is an
element of $\PP_n$ if and only if $\varphi$ is invariant under right
multiplication with the transposition $\tau\in S_n$ swapping $1$ and
$2$. Clearly, this implies again that $\PP_n$ is a left ideal of
$\DD_n$. Furthermore, extensions of $\PP_n$ may now be obtained easily
as in the case of $\DD_n$.  The combinatorial approach extends to
arbitrary finite Coxeter groups.  As a by-product, a short proof of
Solomon's result in the general case is obtained.

As an outer sub-algebra of $\DD$, the peak algebra $\PP$ is also
generated by a series of divided powers and therefore --- as
a general fact --- a homomorphic image of the Hopf algebra $\DD$.
This is shown in Section~\ref{divided} and allows to complete the
picture of duality to Stembridge's peak algebra. In Section~\ref{c},
the commutative counterpart of the peak algebra $\PP_n$ is computed,
which is the dual way of stating \cite[Theorem~3.8]{stembridge97} and
serves as a helpful tool for what follows.

The investigations in $\PP_n$ (as a sub-algebra of $KS_n$) start in
the important Section~\ref{PLI}, where Lie idempotents in $\PP_n$ are
constructed and peak counterparts of the classical Lie idempotents are
displayed.  As a consequence, it is shown that the primitive Lie
algebra of $\PP$ is free, and that $\PP$ is its enveloping algebra.
Some applications related to the Dynkin peak operator and the Klyachko
peak Lie idempotent are presented.  Furthermore, considering peak Lie
idempotents allows to derive an algebraic characterization of $\PP_n$
in Section~\ref{alg-approach}, which rests upon the remarkable
corresponding result on $\DD_n$ given in \cite{garsia-reutenauer89}:
an element $\varphi\in\DD_n$ is contained in $\PP_n$ if and only if
--- via Polya action --- $\varphi$ annihilates the product $P_1\cdots
P_k$ of the homogeneous Lie monomials $P_1,\ldots,P_k$ in the free
associative algebra over a set $X$ whenever $P_1$ is homogeneous of
\emph{even} degree.

An Eulerian-type sub-algebra of $\PP_n$, linearly generated by the
sums of all permutations with a given \emph{number} of peaks, is
introduced and studied in Section~\ref{euler}. Once more, there is a
perfect analogy to the descent case.  A related result is due to Doyle
and Rockmore (\cite{doyle-rockmore02}).

In Section~\ref{radical}, basic information on the structure of
$\PP_n$ can be derived easily from the corresponding results on
$\DD_n$ obtained in \cite{blelau96,blelau02}, since $\PP_n$ turns out
to be a \emph{direct summand} of the regular $\DD_n$-left module.
This includes a description of the Jacobson radical of $\PP_n$ and its
nil-potency index, of the principal indecomposable and the irreducible
$\PP_n$-modules, and of the Cartan invariants.  Finally, there is a
conjecture on the descending Loewy series of $\PP_n$, which is
partially proven.
 
Throughout, $\dim\, M$ denotes the $K$-dimension of a $K$-vector
space $M$, and products $\pi\sigma$ of permutations $\pi,\sigma\in
S_n$ are to be read from left to right: first $\pi$, then $\sigma$.

%
%
%
%
%
%
%
%
%
%
%
%
%
%
%
%
%
%
%
%
%
%
%
%
%

\section{The Solomon descent algebra, quasi-symmetric functions and duality}
\label{descent}

In this section, some notations are fixed and enough of the basic
results on the Solomon descent algebra is reviewed to start our
investigations. A short recall of the duality between this algebra and
the algebra of quasi-symmetric functions follows, allowing to deduce
that the direct sum $\PP$ of all peak algebras $\PP_n$ is the dual of
the algebra of peak functions introduced by Stembridge in
\cite{stembridge97}.

Let $\N$ (respectively, $\N_0$) be the set of all positive (respectively,
nonnegative) integers and put
$$
\haken{n} := \setofall{k\in\N}{k\le n} 
\quad\mbox{ and }\quad
\haken{n}_0 := \setofall{k\in\N_0}{k\le n}
$$
for all integers $n$.  A second basis of $\DD_n$ is constituted by
the elements
$$
\Xi^D := \sum_{E\subseteq D} \Delta^E
\qquad(D\subseteq\haken{n-1}).
$$
By inclusion/exclusion, 
\begin{equation}
  \label{eq:1}
\Delta^D
=
\sum_{E\subseteq D} (-1)^{|D|-|E|}\, \Xi^E,
\end{equation}
for all $D\subseteq\haken{n-1}$.

Sometimes it is more convenient to use compositions of $n$ instead of subsets
of $\haken{n-1}$ as indices. Let $\N^*$ be a free monoid over the alphabet
$\N$. The concatenation product of $q,r\in\N^*$ is denoted by $q.r$ in order to
avoid confusion with the ordinary product in $\N$.
For any $q=q_1\wldots q_k\in \N^*$, 
$\length(q):=k$ is the \emph{length} and $\summe q:=q_1+\cdots+q_k$ is the
\emph{sum} of $q$.  If $\summe q=n$, $q$ is a \emph{composition} of
$n$ ($q\zerl n$).  The mapping
$$
q\lmt D(q):=\{q_1,q_1+q_2,\ldots,q_1+\cdots+q_{k-1}\}
$$
is a bijection from the set of all compositions of $n$ onto the power set of
$\haken{n-1}$, and $\length(q)=|D(q)|+1$ for all $q\zerl n$.  Set
$$
\Delta^q:=\Delta^{D(q)} 
\quad\mbox{ and }\quad
\Xi^q:=\Xi^{D(q)}.
$$
If $r=r_1\wldots r_l\in\N^*$ such that $\summe\,r=\summe\,q$ and
$D(r)\subseteq D(q)$, then $q$ is called a composition of $r$ ($q\zerl
r$).  Equivalently, there exist $q^{(1)},\ldots,q^{(l)}\in\N^*$ such
that $q=q^{(1)}\wldots q^{(l)}$ and $q^{(i)}\zerl r_i$ for all
$i\in\haken{l}$. Now \kref{eq:1} reads as
\begin{equation}
  \label{eq:2}
\Delta^q
=
\sum_{q\zerl r} (-1)^{\length(q)-\length(r)}\, \Xi^r.
\end{equation}

The direct sum
$$
KS := \bigoplus_{n\ge 0} KS_n
$$
is equipped with an \emph{inner product} $\inn$ inherited from the
ordinary products in $KS_n$, $n\ge 0$, by orthogonal extension:
\begin{equation}
  \label{eq:bonus}
  \alpha\inn\beta
  :=
  \tiptop{\alpha\beta}{\mbox{if $n=m$,}}{0}{\mbox{if $n\neq m$,}}
\end{equation}
for all $n,m\in\N_0$, $\alpha\in KS_n$, $\beta\in KS_m$.
Solomon's result mentioned in the introduction may be restated as follows.

\begin{theorem} \label{sol-result}
  The direct sum $ \DD = \bigoplus_{n\ge 0} \DD_n $ is a sub-algebra
  of $(KS,\inn)$.
\end{theorem}

Malvenuto and Reutenauer introduced the structure of a Hopf algebra
on $KS$ (\cite{malvenuto-reutenauer95}), based on the
convolution or \emph{outer product} $\conv$
and a coproduct $\dank$.  The identity element in
$(KS,\conv)$ is $\emptyset$, the empty permutation in $S_0$.
For our purposes, it is enough to recall the following result (ibid.,
Theorem~3.3):

\begin{theorem} \label{bialg-D}
  $\DD$ is a Hopf sub-algebra of $(KS,\conv,\dank)$. Furthermore,
  \begin{equation}
    \label{eq:3}
    \Xi^r\conv\Xi^q
    =
    \Xi^{r.q},
  \end{equation}
  for all $q,r\in\N^*$, and
  \begin{equation}
    \label{eq:4}
    \Xi^n\dank
    =
    \sum_{k=0}^n \Xi^k\tensor \Xi^{n-k}
  \end{equation}
  for all $n\in\N_0$, where $\Xi^0:=\emptyset$.  In particular,
  $(\DD,\conv)$ is a polynomial ring in the (non-commuting) variables
  $$
  \Xi^n=\sum_{\Des(\pi)=\emptyset} \pi=\id_n\qquad
  (n\in\N),
  $$
  where $\id_n$ denotes the identity in $S_n$.
\end{theorem}

The Hopf algebra $(\DD,\conv,\dank)$ is isomorphic to the Hopf algebra
of noncommutative symmetric functions (\cite{parisI}).  In that
setting a crucial connection between the inner and the outer product
on $\DD$ has been established (\cite[Proposition~5.2]{parisI}). To
state it, note that the inner product $\inn$ on $\DD$ canonically
gives rise to a (component-wise) product on $\DD\tensor\DD$ (which is
also denoted by $\inn$).

\begin{lemma}
  \label{reziprozi}
  For all $\alpha,\beta,\gamma\in\DD$, 
  $$
  (\alpha\conv\beta)\inn\gamma =
  ((\alpha\tensor\beta)\inn(\gamma\dank))\upnk\,,
  $$
  where $\upnk:\DD\tensor\DD\lra\DD$ is the linearization $
  (\varphi\tensor\psi)\upnk = \varphi\conv\psi $ of the outer product
  $\conv$ on $\DD$.
\end{lemma}

\breakitdown

Let $X=\{x_1,x_2,\ldots\}$ be an infinite set of commuting variables
and denote by
$$
M_r
:=
\sum_{i_1<\cdots<i_l} x_{i_1}^{r_1}\cdots x_{i_l}^{r_l}
$$
the \emph{monomial quasi-symmetric function} indexed by 
$r$, for all $r=r_1\wldots r_l\in\N^*$.
The $M_r$'s constitute a linear basis of the algebra $\Qsym$ of quasi-sym\-metric
functions introduced by Gessel (\cite{gessel84}).
A second basis of $\Qsym$ consists of the \emph{fundamental quasi-symmetric
  functions}, defined by
$$
F_q
:=
\sum_{r\zerl q} M_r
\qquad 
(q\in\N^*).
$$
As before, subsets of $\haken{n-1}$ as indices work equally well
(whenever $n$ is clear from the context): For $D\subseteq\haken{n-1}$,
choose $r\zerl n$ such that $D=D(r)$ and set
$$
M_D
:=
M_r 
\quad\mbox{ and }\quad
F_D
:=
F_r\,.
$$
The algebra $\Qsym$ is equipped with two coproducts, the
\emph{outer coproduct} $\gamma$ and the \emph{inner coproduct}
$\gamma'$; the former turns $\Qsym$ into a (graded) Hopf algebra.  For
details, the reader is referred to \cite{malvenuto-reutenauer95}.  In
the same paper, there is the following result which allows to transfer
results on quasi-symmetric functions to the Solomon algebra and vice
versa (ibid., Theorems~3.2 and~3.3).

\begin{theorem}  \label{major-duality}
  The graded dual $\Qsym^*$ of $\Qsym$ with the product inherited from the
  outer coproduct $\gamma$ and the coproduct inherited from the ordinary
  product of $\Qsym$, is isomorphic to the Hopf algebra $(\DD,\conv,\dank)$.
  
  In this isomorphism, the basis $\setofall{F_q^*}{q\in\N^*}$ of $\Qsym^*$ dual
  to the basis $\setofall{F_q}{q\in\N^*}$ of $\Qsym$ corresponds to the
  basis $\setofall{\Delta^q}{q\in\N^*}$ of $\DD$.
  
  Furthermore, the same mapping is an isomorphism of the algebra 
  $\Qsym^*$ with the product inherited from the inner coproduct
  $\gamma'$ of $\Qsym$, onto $(\DD,\inn)$.
\end{theorem}
 
The latter part of the preceding theorem is due to Gessel.
Define a pairing between $\DD$ and $\Qsym$ by
$$
\pairing{\Delta^q}{F_s}
:=
\tiptop{1}{\mbox{if $q=s$,}}{0}{\mbox{if $q\neq s$,}}
$$
then the above result leads to the following reciprocity laws:
\begin{eqnarray*}
\pairing{\varphi\inn\psi}{f}
& = &
\pairing{\varphi\tensor\psi}{\gamma'(f)}\\[1mm]
\pairing{\varphi\conv\psi}{f}
& = &
\pairing{\varphi\tensor\psi}{\gamma(f)}\\[1mm]
\pairing{\varphi}{fg}
& = &
\pairing{\varphi\dank}{f\tensor g}
\end{eqnarray*}
for all $\varphi,\psi\in\DD$, $f,g\in\Qsym$.  Furthermore,
Theorem~\ref{major-duality} implies that $\dank$ is also a
homomorphism for the inner products on $\DD$ and $\DD\tensor\DD$.

\breakitdown

A brief argument follows which shows that $\PP$ is the dual of the
sub-algebra $\PPP$ of $\Qsym$ introduced and studied by Stembridge in
\cite{stembridge97}.

Set $D+z:=\setofall{d+z}{d\in D}$ for all integers $z$ and
$D\subseteq\N_0$, then
\begin{equation}
  \label{eq:5}
  \Peak(\pi)
  =
  \Des(\pi)\cap \Big((\haken{n-1}\setminus\Des(\pi))+1\Big),
\end{equation}
for each $\pi\in S_n$.  In particular, $P\subseteq\haken{n}$ is a peak
set of a permutation $\pi\in S_n$ if and only if
$P\subseteq\{2,\ldots,n-1\}$ and $P\cap (P-1)=\emptyset$.  In this
case, $P$ is called a \emph{peak set in $\haken{n}$}.

Each homogeneous component $\PPP_n$ of $\PPP$ has a linear basis constituted by
the elements
$$
K_P
:=
\sum_{\substack{D\subseteq\haken{n-1}\\[2pt] P\subseteq D\scsctriangle(D+1)}}
F_D\,,
$$
indexed by peak sets $P$ in $\haken{n}$, where $D\sctriangle
(D+1)=D\backslash(D+1)\cup (D+1)\backslash D$ denotes the symmetric
difference (\cite[Theorem~3.1(a) and Proposition~3.5]{stembridge97}).
Furthermore, if $P(D):=D\cap((\haken{n-1}\backslash D)+1)$ for all
$D\subseteq\haken{n-1}$, then the linear mapping
$\theta:\Qsym\to\PPP$, defined by
$$
F_D\lmt K_{P(D)}\,,
$$
is an epimorphism of algebras 
(\cite[Theorem~3.1(c) and its proof]{stembridge97}).
Now, defining
$$
\PP^\perp
:=
\setofall{f\in\Qsym}{\pairing{\varphi}{f}=0\mbox{ for all }\varphi\in\PP},
$$
we are ready to state and prove the following duality, however, for the time
being, with the restriction to the coproduct in $\DD$ and the product in
$\Qsym$.

\begin{mtheorem} \label{duality-1}
  $\PP^\perp=\ker\,\theta$.

  In particular, $\PP$ is a sub-coalgebra of $(\DD,\dank)$, dual to the
  Stembridge peak algebra $\PPP$, and
  $$
  [\Pi^P,K_Q]
  :=
  \tiptop{1}{\mbox{if $P=Q$,}}{0}{\mbox{otherwise,}}
  \qquad
  \mbox{($P$, $Q$ peak sets in $\haken{n}$, $n\in\N_0$)}
  $$
  defines a pairing between $\PP$ and $\PPP$ such that
  $$
  [\varphi,fg]=[\varphi\dank,f\tensor g]
  $$
  for all $\varphi\in\PP$, $f,g\in\PPP$.
\end{mtheorem}

\begin{proof}
  For any element $f=\sum_D a_D F_D\in\Qsym_n$, $ f\perp \PP $ if and only if,
  for all peak sets $P$ in $\haken{n}$,
  $$
  0
  =
  \pairing{\Pi^P}{f}
  =
  \sum_D a_D \sum_{P(E)=P} \pairing{\Delta^E}{F_D}
  =
  \sum_{P(D)=P} a_D\,.
  $$
  But this is equivalent to
  $$
  f\theta
  =
  \sum_D a_D F_D\theta
  =
  \sum_P \Big(\sum_{P(D)=P} a_D\Big) K_P
  =
  0,
  $$
  that is, $f\in\ker\,\theta$.  This proves the first claim.  As a
  consequence, $\PP^\perp$ is an ideal of $\Qsym$, hence $\PP$ a
  sub-coalgebra of $\DD$, by Theorem~\ref{major-duality}.
  Furthermore, as a general fact, $\PP$ is the graded dual of
  $\Qsym/\PP^\perp$, with pairing
  $$
  [\varphi,\ol{f}\,]:=\pairing{\varphi}{f}
  $$
  for all $\varphi\in\PP$, $f\in\Qsym$, where
  $\ol{f}:=f+\PP^\perp$.  The observation that the mapping
  $\ol{\theta}:\ol{F}_D\mapsto K_{P(D)}$ is an isomorphism of algebras
  from $\Qsym/\ker\,\theta=\Qsym/\PP^\perp$ onto $\PPP$, completes the
  proof.
\end{proof}

The situation dual to the one in the preceding theorem --- displaying
$\PP$ as a factor algebra of $(\DD,\conv)$ and dual to $\PPP$ as a
sub-coalgebra of $(\Qsym,\gamma)$ --- is analyzed in Section~\ref{divided}.
However, the duality of $\PP$ and $\PPP$ is not needed in what
follows.  Main Theorem~\ref{duality-1} rather serves as a link to the
theory of quasi-symmetric functions for the reader more familiar with
that setting.

%
%
%
%
%
%
%
%
%
%
%
%
%
%
%
%
%
%
%
%
%
%
%
%
%

\section{The peak algebra $\PP$}
\label{peak}

The investigations in the peak algebra $\PP = \bigoplus_{n\ge 0}
\PP_n$ now begin with establishing a triple algebraic structure on
$\PP$ as in the case of the descent algebra.  $\PP$ is a sub-bialgebra
of $(\DD,\conv,\dank)$ and --- as a consequence of this --- a left
ideal of $(\DD,\inn)$.
 
The number of peak sets in $\haken{n}$ (and thus the dimension of
$\PP_n$) is readily seen to be the Fibonacci number $f_n$, defined by
$f_1=f_2=1$ and
$$
f_n=f_{n-1}+f_{n-2}
$$
for all $n\ge 3$.  This number equals the number of \emph{odd
  compositions} $q$ of $n$ ($q\zerlodd n$), where $q=q_1\wldots
q_k\in\N^*$ is called \emph{odd} if $q_i$ is odd, for all
$i\in\haken{k}$. To obtain a correspondence between peak sets in
$\haken{n}$ and odd compositions $q$ of $n$, choose
$m_1,\ldots,m_k\in\N_0$ such that $q_i=2m_i+1$, for all
$i\in\haken{k}$, and set
$$
\tilde{q} := 2^{.m_1}.1.2^{.m_2}.1\wldots 2^{.m_k}.1\zerl q,
$$
where $d^{.m}$ denotes the $m$-th power of $d$ in $\N^*$, for all
$d\in\N$, $m\in\N_0$.  Observing that
\begin{equation}
  \label{eq:7}
  D(q)=D(\tilde{q})\cap (D(\tilde{q})+1)
\end{equation}
we claim the following.

\begin{proposition} \label{peak-comp}
  The mapping
  $$
  q \lmt P(q) := (\haken{n-1}\setminus D(\tilde{q}))+1
  $$
  is a one-one correspondence between the odd compositions of $n$
  and the peak sets in $\haken{n}$, and
  \begin{equation}
    \label{eq:8}
    D(q)=\haken{n-1}\setminus((P(q)-1)\cup P(q)),    
  \end{equation}
  for all $q\zerlodd n$.
\end{proposition}

\begin{proof}
  Equation \kref{eq:8} is immediate from equation \kref{eq:7}.
  
  Furthermore, $1\notin P(q)$ is obvious, while $n-1\in D(\tilde{q})$
  implies that $n\notin P(q)$, which shows $P(q)\subseteq\{2,\ldots,n-1\}$.
  From $\tilde{q}_j\le 2$ for all $j$, it follows that $P(q)\cap
  (P(q)-1)=\emptyset$, hence $P(q)$ is a peak set in $\haken{n}$.
  
  If $q,r\zerlodd n$, then $P(q)=P(r)$ implies $D(q)=D(r)$, by
  \kref{eq:8}, thus $q=r$.
  
  The proof is complete, since both sets in question have cardinality
  $f_n$.
\end{proof}

\breakitdown

To analyze outer products of the elements of $\PP$ (based on formula
\kref{eq:3}), consider the basis of $\PP_n$ constituted by the elements
\begin{equation}
  \label{triangle}
  \Gamma^P
  :=
  \sum_Q \Pi^Q
  \qquad\mbox{($P$ a peak set in $\haken{n}$)},  
\end{equation}
where the sum is taken over all peak sets $Q$ in $\haken{n}$ containing $P$.
By inclusion/exclusion,
\begin{equation}
  \label{triangle-back}
  \Pi^P
  =
  \sum_Q (-1)^{|Q|-|P|}\Gamma^Q
\end{equation}
with the same range of the sum.  Proposition~\ref{peak-comp} allows to take odd
compositions instead of peak sets as indices:
$$
\Gamma^q := \Gamma^{P(q)},
$$
for all $q\zerlodd n$.

\begin{proposition} \label{gamma-xi}
  The elements $\Gamma^r$, $r\zerlodd n$, constitute a linear basis
  of $\PP_n$.  Furthermore,
  $$
  \Gamma^q = (-1)^{(n-\length(q))/2} \sum_{\tilde{q}\zerl s\zerl q}
  (-1)^{\length(s)-\length(q)}\; \Xi^s
  $$
  for all $q\zerlodd n$.
\end{proposition}

\begin{proof}
  The basis property follows from Proposition~\ref{peak-comp}.
  
  Let $q\zerlodd n$ and $P:=P(q)$, then \kref{eq:5} implies that
  $P\subseteq \Peak(\pi)$ if and only if $P\subseteq
  \Des(\pi)\subseteq \haken{n-1}\setminus(P-1)$, for all $\pi\in S_n$.
  Now combine \kref{eq:1} and \kref{eq:8} with standard set theoretic
  arguments to obtain
  \begin{eqnarray*}
    \Gamma^P
    & = &
    \sum_{%
      \substack{Q\;\mbox{\footnotesize peak set in}\;\haken{n}\\[2pt] 
        P\,\subseteq\, Q}}
    \Pi^Q\\[2mm]
    & = &
    \sum_{P\,\subseteq\, D\,\subseteq\, \haken{n-1}\setminus(P-1)} 
    \Delta^D \\[2mm]
    & = &
    \sum_{P\,\subseteq\, D\,\subseteq\, \haken{n-1}\setminus(P-1)} 
    \quad
    \sum_{E\,\subseteq\, D}
    (-1)^{|D|-|E|}\,\; \Xi^E \\[2mm]
    & = &
    \sum_{E\,\subseteq\, \haken{n-1}\backslash (P-1)}
    \quad
    \sum_{E\cup P\,\subseteq\, D\,\subseteq\, \haken{n-1}\setminus(P-1)} 
    (-1)^{|D|-|E|}\,\; \Xi^E \\[2mm]
    & = &
    \sum_{E\,\subseteq\, \haken{n-1}\backslash (P-1)}
    (-1)^{|E\cup P|-|E|}\,\; \Xi^E
    \sum_{\tilde{D}\,\subseteq\, (\haken{n-1}\setminus(P-1))\setminus(E\cup P)}
    (-1)^{|\tilde{D}|}  \\[2mm]
    & = &
    \sum_{%
      \haken{n-1}\setminus ((P-1)\cup P)\,
      \subseteq\, 
      E\,
      \subseteq\, 
      \haken{n-1}\setminus (P-1)}
    (-1)^{|E\cup P|-|E|}\;\, \Xi^E \\[2mm]
    & = &
    \sum_{D(q)\,\subseteq\, E\,\subseteq\, D(\tilde{q})}
    (-1)^{|E\cup P|-|E|}\,\; \Xi^E. 
  \end{eqnarray*}
  But, for each $E\subseteq\haken{n-1}$ such that
  $D(q)\subseteq E\subseteq D(\tilde{q})$,
  $$
  E\setminus P
  \subseteq
  D(\tilde{q})\setminus P
  =
  D(\tilde{q})\cap (D(\tilde{q})+1)
  \subseteq
  D(q),
  $$
  hence $|E\cup P|=|D(q)\cup P|=|D(q)|+|P|$.
  It follows that
  \begin{eqnarray*}
  \Gamma^P
  & = &
  (-1)^{|P|}
    \sum_{D(q)\,\subseteq\, E\,\subseteq\, D(\tilde{q})}
    (-1)^{|E|-|D(q)|}\,\; \Xi^E\\[1mm]
  & = &
  (-1)^{(n-\length(q))/2}
  \sum_{\tilde{q}\zerl s\zerl q}
    (-1)^{\length(s)-\length(q)}\,\; \Xi^s        
  \end{eqnarray*}
  as asserted.
\end{proof}

A direct consequence of Proposition~\ref{gamma-xi} is

\begin{proposition} \label{out-prod-P}
  $(\PP,\conv)$ is a sub-algebra of $(\DD,\conv)$.  Furthermore,
  $$
  \Gamma^r\conv\Gamma^q = \Gamma^{r.q}
  $$
  for all $q, r\in\N^*$ odd.
\end{proposition}

\begin{proof}
  For all $q\zerlodd  n$, $r\zerlodd  m$, by Proposition~\ref{gamma-xi} and
  \kref{eq:3},
  \begin{eqnarray*}
  \Gamma^r\conv\Gamma^q
  & = &
  (-1)^{(m+n-\length(r)-\length(q))/2}
  \sum_{\tilde{r}\zerl s\zerl r}
  \sum_{\tilde{q}\zerl t\zerl q}
  (-1)^{\length(s)-\length(r)+\length(t)-\length(q)}\, 
  \Xi^s\conv\Xi^t\\
  & = &
  (-1)^{(m+n-\length(r.q))/2}
  \sum_{\tilde{r}.\tilde{q}\zerl u\zerl r.q} 
  (-1)^{\length(u)-\length(r.q)}\, \Xi^u\\
  & = &
  \Gamma^{r.q},
  \end{eqnarray*}
  since $\tilde{r}.\tilde{q}=\widetilde{r.q}$.
\end{proof}

The preceding result is the dual way of stating
\cite[Theorem~2.2]{bmsw02} (see also Corollary~\ref{duality-2}).  It
seems to be worth mentioning that $(\PP,\conv)$ is thus a polynomial
ring in the (non-commuting) variables
$\Gamma^n=\Pi^{\{2,4,\ldots,n-1\}}$, $n\in\N$ odd, each of which is
the sum of all permutations $\pi\in S_n$ with the maximal peak set
$\{2,4,\ldots,n-1\}$. These permutations are sometimes called
\emph{alternating}, since
$$
1\pi<2\pi>3\pi<4\pi>\cdots<(n-1)\pi>n\pi.
$$

\breakitdown

To analyze co-products of elements of $\PP_n$, it is more convenient
to consider the minimal peak set $\emptyset$: Define
$$
\R^n = 2\Pi^\emptyset \in \PP_n,
$$
for all $n\in\N$. This is (up to the factor $2$) the sum of the \emph{valley
  permutations} $\pi$ in $S_n$, defined by the property
$$
1\pi>2\pi>\cdots>k\pi<(k+1)\pi<\cdots<n\pi,
$$
where $k:=1\pi^{-1}$.
The factor $2$ and the multiplication rule
\begin{equation}
  \label{vert-hor}
  \Delta^{1^{.k}}\conv\Delta^m = \Delta^{1^{.k}.m}+\Delta^{1^{.(k-1)}.(m+1)}
\end{equation}
for all $k,m\in\N$ (\cite[(52)]{parisI}) allow to rewrite the
definition to
\begin{equation}
  \label{re-xi-tilde}
  \R^n
  =
  2\sum_{k=0}^n \Delta^{1^{.k}.(n-k)}
  =
  \sum_{k=0}^n \Delta^{1^{.k}}\conv\Delta^{n-k},
\end{equation}
where $\Delta^0:=\emptyset$.  For all
$q=q_1\wldots q_k\in\N^*$, set
$$
\R^q := \R^{q_1} \conv \cdots \conv \R^{q_k},
$$
then $\R^q\in\PP$, by Proposition~\ref{out-prod-P}, and these
elements may be expressed in terms of the $\Pi$-basis and the
$\Gamma$-basis of $\PP$, and the $\Xi$-basis of $\DD$, as follows.

For all $q=q_1\wldots q_k\zerl n$, let $q^\dagger := q_k$ denote the
final letter of $q$. More generally, if $q\zerl r=r_1\wldots r_l$, say
$q=q^{(1)}\wldots q^{(l)}$ such that $q^{(i)}\zerl r_i$ for all
$i\in\haken{l}$, let
$$
F_r(q) := \prod_{i=1}^l (q^{(i)})^\dagger
$$
be the product of the \emph{$r$-finals of $q$}.

\begin{proposition} \label{transition-tilde-xi}
  Let $n\in\N$, $r\zerl n$ and put $D:=D(r)\cup (D(r)+1)$, then
  \begin{eqnarray*}
    \R^r
    & = &
    2^{\length(r)}
    \sum_{%
      \substack{Q\;\mbox{\footnotesize peak set in}\;\haken{n}\\[3pt]
        Q\subseteq D}}
    \Pi^Q\\[1mm]
    & = &
    2^{\length(r)}
    \sum_{%
      \substack{Q\;\mbox{\footnotesize peak set in}\;\haken{n}\\[3pt]
        Q\subseteq \haken{n-1}\backslash D}}
    (-1)^{|Q|}\; \Gamma^Q\\[1mm]
    & = &
    2^{\length(r)}
    \sum_{\substack{q\zerl r\\[2pt] q\;\mbox{\footnotesize odd}}}
    (-1)^{(n-\length(q))/2}\; \Gamma^q\\
    & = &
    2^{\length(r)}
    \sum_{\substack{q\zerl r\\[2pt] F_r(q)\;\mbox{\footnotesize odd}}}
    (-1)^{n-\length(q)}\; \Xi^q.
  \end{eqnarray*}
  In particular,
  \begin{eqnarray*}
    \R^n
    & = &
    2
    \sum_{q\zerlodd\, n} (-1)^{(n-\length(q))/2}\; \Gamma^q
      = 
    2
    \sum_{\substack{q\zerl n\\[2pt] q^\dagger\;\mbox{\tiny odd}}} 
    (-1)^{n-\length(q)}\; \Xi^q.
  \end{eqnarray*}  
\end{proposition}

\begin{proof}
  To begin with, consider the case $r=n$, then
  \kref{triangle-back} implies
  $\R^n=2\Pi^\emptyset=2\sum_{Q\subseteq \haken{n-1}} (-1)^{|Q|} \Gamma^Q$,
  hence $\R^n=2\sum_{q\zerlodd  n} (-1)^{(n-\length(q))/2} \Gamma^q$, by
  Proposition~\ref{peak-comp}.  
  Furthermore, by \kref{re-xi-tilde} and \kref{eq:2},
  \begin{eqnarray*}
    \R^n
    & = &
    \sum_{k=0}^n \Delta^{1^{.k}}\conv\Xi^{n-k}\\
    & = &
    \sum_{k=0}^{n-1} \sum_{r\zerl k} (-1)^{k-\length(r)}\;\Xi^{r.(n-k)}
    +
    \sum_{s\zerl n} (-1)^{n-\length(s)}\;\Xi^s\\
    & = &
    \sum_{s\zerl n} 
    ((-1)^{n-s^\dagger-(\length(s)-1)}+(-1)^{n-\length(s)})\;\Xi^s\\
    & = &
    \sum_{s\zerl n} 
    (-1)^{n-\length(s)}((-1)^{s^\dagger+1}+1)\;\Xi^s\\
    & = &
    2\sum_{\substack{s\zerl n\\[2pt] s^\dagger\;\mbox{\tiny odd}}} 
    (-1)^{n-\length(s)}\;\Xi^s
  \end{eqnarray*}
  as asserted. Now let $r=r_1\wldots r_l$, then the latter identity
  implies
  \begin{eqnarray*}
    \R^r
    & = &
    \R^{r_1}\conv\cdots\conv\R^{r_l}\\
    & = &
    \Big(
    2\sum_{\substack{s\zerl r_1\\[2pt] s^\dagger\;\mbox{\tiny odd}}}
    (-1)^{r_1-\length(s)}\;\Xi^{s}
    \Big)
    \conv
    \cdots
    \conv
    \Big(
    2\sum_{\substack{s\zerl r_l\\[2pt] s^\dagger\;\mbox{\tiny odd}}} 
    (-1)^{r_l-\length(s)}\;\Xi^{s}
    \Big)\\
    & = &
    2^l
    \sum_{\substack{s\zerl r\\[2pt] F_r(s)\;\mbox{\tiny odd}}} 
    (-1)^{n-\length(s)}\;\Xi^s.
  \end{eqnarray*}
  Similarly,
  \begin{eqnarray*}
    \R^r
    & = &
    \Big(
    2\sum_{q\zerlodd  r_1} 
    (-1)^{(r_1-\length(q))/2}\;\Gamma^q
    \Big)
    \conv
    \cdots
    \conv
    \Big(
    2\sum_{q\zerlodd  r_l} 
    (-1)^{(r_l-\length(q))/2}\;\Gamma^q
    \Big)\\
    & = &
    2^l
    \sum_{\substack{q\zerl r\\[2pt] q\;\mbox{\tiny odd}}}
    (-1)^{(n-\length(q))/2}\;\Gamma^q.
  \end{eqnarray*}
  But, for each $q\zerlodd n$, the condition $q\zerl r$ is equivalent to
  $$
  D(r)\subseteq D(q)=\haken{n-1}\backslash((P(q)-1)\cup P(q)),
  $$
  by \kref{eq:8}, and this is another way of stating $P(q)\subseteq
  \haken{n-1}\backslash D$.  Finally,
  $$
  \sum_{Q\subseteq \haken{n-1}\backslash D} (-1)^{|Q|}\Gamma^Q =
  \sum_{Q\subseteq D} \Pi^Q
  $$
  easily follows from \kref{triangle}.
\end{proof}

We are now in a position to state and prove the peak counterpart of
Theorem~\ref{bialg-D}.

\begin{mtheorem} \label{bialg-P}
  $\PP$ is a Hopf sub-algebra of $(\DD,\conv,\dank)$.  Furthermore,
  the elements $\R^q$, $q\in\N^*$ odd, constitute a linear basis of
  $\PP$, and
  $$
  \R^r\conv\R^q = \R^{r.q},
  $$
  for all $q,r\in\N^*$, while
  $$
  \R^n\dank = \sum_{k=0}^n \R^k\tensor \R^{n-k},
  $$
  for all $n\in\N$, where $\R^0:=\emptyset$.  
  In particular, $(\PP,\conv)$ is a polynomial
  ring in the (non-commuting) variables
  $$
  \textstyle\frac{1}{2}\,\R^n
  =
  \displaystyle\sum_{\Peak(\pi)=\emptyset} \pi\qquad (n\in\N\mbox{ odd}).
  $$
\end{mtheorem}

\begin{proof}
  Expressing the elements $\R^q$, $q\zerlodd n$, in the $\Xi$-basis of
  $\DD_n$ yields linear independency, by
  Proposition~\ref{transition-tilde-xi} and triangularity.  These
  elements thus constitute a linear basis of $\PP_n$ as asserted.  The
  multiplication rule is $\R^r\conv\R^q = \R^{r.q}$, by definition,
  while the co-multiplication rule follows from standard arguments
  based on \kref{eq:4}, \kref{re-xi-tilde}, and
  $$
  \Delta^{1^{.n}}\dank 
  = 
  \sum_{k=0}^n \Delta^{1^{.k}}\tensor\Delta^{1^{.(n-k)}}
  $$
  for all $n\in\N$ (see \cite[3.8]{parisI}).
\end{proof}

A sub-algebra $\AA$ of $(KS,\conv)$ is \emph{graded by intersection}
if $\AA=\bigoplus_{n\ge 0} (\AA\cap KS_n)$.

Main Theorem~\ref{bialg-P} combined with the following consequence of
Lemma~\ref{reziprozi} yields a \emph{first proof of Main
  Theorem~\ref{left-ideal}}.

\begin{proposition}\label{leftideal}
  Let $\AA$ be a sub-bialgebra of $(\DD,\conv,\dank)$ and assume that
  $\AA$ is graded by intersection, then $\AA$ is a left ideal of $\DD$
  with respect to the inner product.
\end{proposition}

\begin{proof}
  Let $q\in\N^*$ and proceed by induction on $\length(q)$ to prove
  that $\Xi^q\inn\AA\subseteq \AA$.  For $\length(q)=1$, this is
  another way of stating that $\AA$ is graded by intersection.  Let
  $n\in\N$, then for all $\alpha\in\AA$, there exist
  $\alpha^{(1)}_i,\alpha^{(2)}_i\in\AA$ ($i\in I$) such that
  $\alpha\dank=\sum_{i\in I} \alpha^{(1)}_i\tensor \alpha^{(2)}_i$,
  hence
  \begin{eqnarray*}
    \Xi^{q.n}\inn\alpha
    & = &
    ((\Xi^q\tensor\Xi^n)\inn(\alpha\dank))\upnk
    = 
    \sum_{i\in I} (\Xi^q\inn\alpha^{(1)}_i)\conv \alpha^{(2)}_i
    \in 
    \AA,
  \end{eqnarray*}
  by Lemma~\ref{reziprozi} and induction.
\end{proof}

%
%
%
%
%
%
%
%
%
%
%
%
%
%
%
%
%
%
%
%
%
%
%
%
%

\section{A Combinatorial Approach} \label{comb-approach}

There is a simple combinatorial explanation of Main
Theorems~\ref{left-ideal} and~\ref{bialg-P}, which will be given in
this section.
  
Blessenohl and Laue introduced a graph structure on $S_n$ to describe
the descent classes $\setofall{\pi\in S_n}{\Des(\pi)=D}$,
$D\subseteq\haken{n-1}\,$, and to derive a combinatorial proof of
Theorem~\ref{sol-result} (\cite[Section~4]{blelau93a}).  A slight
extension of this graph structure yields a combinatorial
characterization of the elements of $\PP_n$ and, as a consequence, a
two-line proof of Main Theorem~\ref{left-ideal}. Furthermore, peak
analogs of certain extensions of the algebra $\DD_n$ may easily be
obtained (see Remark~\ref{extensions}).
  
This combinatorial concept extends to arbitrary finite Coxeter groups
(instead of $S_n$) and allows a short proof of Solomon's result in the
general case, as is briefly explained at the end of this section.
Arising from this, there is the notion of a peak-like sub-algebra of
the group algebra of a finite Coxeter group.
  
Throughout, $\tau_{n,i}$ denotes the transposition swapping $i$ and  
$i+1$ in $S_n$, for all $n\in\N$ and $i\in\haken{n-1}$, and  
$\tau_{1,1}$ denotes the identity in $S_1$.  All what follows is based  
on the observation that  
\begin{equation}  
  \label{1-2-inv}  
  \Peak(\pi)  
  =  
  \Peak(\pi\tau_{n,1})
\end{equation}  
for all $\pi\in S_n$, that is, swapping $1$ and $2$ in the image line  
of $\pi$ does not change the peak set of $\pi$.  Indeed, if  
$\{1\pi^{-1},2\pi^{-1}\}=\{k,l\}$ such that $k<l$, then either  
$l-k>1$, hence even $\Des(\pi)=\Des(\pi\tau_{n,1})$, or $l=k+1$. In  
the latter case, $\Des(\pi)$ and $\Des(\pi\tau_{n,1})$ differ exactly  
by the element $k$, which also implies \kref{1-2-inv}, since  
$k-1\in\Des(\pi)\cup\{0\}$ and $k+1\notin\Des(\pi)$.  
  
As a consequence of \kref{1-2-inv}, $\Pi^P\tau_{n,1} = \Pi^P$ for
all peak sets $P$ in $\haken{n}$, which is one part of
  
\begin{mtheorem} \label{comb-char}  
  Let $\varphi\in\DD_n$, then $\varphi\in\PP_n$ if and only if
  $\varphi\tau_{n,1}=\varphi$.
\end{mtheorem}

This characterization implies, in particular, that $\PP_n$ is a left
ideal of $\DD_n$, that is, \emph{a second proof of Main
  Theorem~\ref{left-ideal}}. For, if $\alpha\in\DD_n$ and
$\varphi\in\PP_n$, then $\alpha\varphi\in\DD_n$, by
Theorem~\ref{sol-result}, and
$(\alpha\varphi)\tau_{n,1}=\alpha(\varphi\tau_{n,1})=\alpha\varphi$,
hence $\alpha\varphi\in\PP_n$.
  
\begin{remark} \label{extensions}  
  Main Theorem~\ref{comb-char} may be restated as
  $$  
  \PP_n = \DD_n\cap \Big(KS_n(\id_n+\tau_{n,1})\Big) = \DD_n\cap  
  \Big(\DD_n(\id_n+\tau_{n,1})\Big).  
  $$  
  The direct sum  
  $  
  KS^{12} := KS_0\oplus\bigoplus_{n\ge 1} KS_n(1+\tau_{n,1})  
  $
  is a sub-algebra of $(KS,\conv)$ and a left co-ideal of
  $(KS,\dank)$, that is, $KS^{12}\dank\subseteq KS^{12}\tensor KS$.
  This is readily seen from the combinatorial description of the
  outer product $\conv$ and the coproduct $\dank$ given in
  \cite[pp. 977--978]{malvenuto-reutenauer95}. It follows that
  $$  
  \PP = \DD\cap KS^{12}  
  $$
  is a convolution sub-algebra and a left co-ideal of $\DD$. Thus,
  since $(\DD,\dank)$ is co-commu\-ta\-tive, $\PP$ is a sub-bialgebra
  of $(\DD,\conv,\dank)$, which essentially recovers Main
  Theorem~\ref{bialg-P}.
  
  Clearly, peak analogs of the \emph{extensions} of $\DD_n$ introduced
  in \cite{loday-ronco98,patrasreut01,schocker-lia} may be obtained
  along the same lines.  For instance, in the above argument, consider
  the sub-algebra $\LL$ (respectively, $\hat{\LL}$) of $(KS,\conv)$
  generated by all Lie elements (respectively, Lie idempotents) in
  $KS_n$, $n\in\N_0$, instead of $\DD$ (see
  \cite{patrasreut01,schocker-lia}, or Section~\ref{PLI}, for the
  necessary definitions). $\LL$ and $\hat{\LL}$ are inner sub-algebras
  of $KS$ (containing $\DD$), and co-commutative Hopf sub-algebras of
  $(KS,\conv,\dank)$ (\cite{patrasreut01,schocker-lia}), hence so are
  $\LL^{12}:=\LL\cap KS^{12}$ and $\hat{\LL}^{12}:=\hat{\LL}\cap
  KS^{12}$; and
  $$  
  \PP 
  = 
  \DD\cap KS^{12}  
  \subseteq
  \hat{\LL}^{12}
  \subseteq
  \LL^{12}.
  $$
  In particular, there is the chain $\PP_n\subseteq
  \hat{\LL}^{12}_n \subseteq \LL^{12}_n$ of sub-algebras of the group
  algebra $KS_n$.
  A similar argument applies to the Hopf algebra of the planar binary
  trees introduced by Loday and Ronco (\cite{loday-ronco98}).
  
  We refrain from a more detailed analysis of the properties of these
  extensions of $\PP_n$.
\end{remark}

The proof of the remaining part of Main Theorem~\ref{comb-char} requires a  
little more work and is based on a \emph{local} description of the  
peak classes  
$$  
\setofall{\pi\in S_n}{\Peak(\pi)=P}  
\qquad  
\mbox{($P$ a peak set in $\haken{n}$).}  
$$

Let $k,l\in\haken{n}$ and $\pi,\sigma\in S_n$, then $k$ is called a  
\emph{neighbor of $l$ in $\pi$} if $|k\pi^{-1}-l\pi^{-1}|=1$, that is, if $k$  
is a neighbor of $l$ in the image line $(1\pi).(2\pi)\wldots(n\pi)$ of  
$\pi$.  
  
Following Blessenohl and Laue, $\sigma$ is called a \emph{descent
  neighbor} of $\pi$ if there is an index $i\in\haken{n-1}$ such that
$\sigma=\pi\tau_{n,i}$ and $i$ is not a neighbor of $i+1$ in $\pi$. In
this case, write $\sigma\descn\pi$.  Denote by $\descr$ the smallest
equivalence relation on $S_n$ containing $\descn$ to recall
\cite[Proposition~4.2]{blelau93a}:
  
\begin{lemma} \label{dbhl-desc}  
  $\Des(\pi)=\Des(\sigma)$ if and only if $\pi\descr\sigma$, for all
  $\pi,\sigma\in S_n$.
\end{lemma}  
  
The identity \kref{1-2-inv} suggests the following analog of
Lemma~\ref{dbhl-desc} for peak sets.  $\sigma$ is called a \emph{peak
  neighbor} of $\pi$ if $\sigma$ is a descent neighbor of $\pi$ or
$\sigma=\pi\tau_{n,1}$. In this case, write $\sigma\peakn\pi$.
Furthermore, denote by $\peakr$ the smallest equivalence relation on
$S_n$ containing $\peakn$.
  
For instance, $\underline{3}\,1\,6\,\underline{2}\,4\,5$ is a descent  
neighbor (hence also a peak neighbor) of $2\,1\,6\,3\,4\,5$, while  
$\underline{1}\,\underline{2}\,6\,3\,4\,5$ is a peak neighbor, but not  
a descent neighbor of $2\,1\,6\,3\,4\,5$.

\begin{proposition} \label{1-2}  
  Let $\pi\in S_n$ and $i\in\haken{n-1}$ such that $i$, $i+1$ are  
  neighbors in $\pi$ and each neighbor $\neq i,i+1$ of $i$ and $i+1$  
  in $\pi$ is larger than $i+1$, then  
  $$  
  \pi \peakr \pi\tau_{n,i}.  
  $$  
\end{proposition}  
  
\begin{proof}  
  By induction on $i$.     
  If $i=1$, then $\pi\peakn\pi\tau_{n,1}$.  
    
  Let $i>1$, then $i-1$ is not a neighbor of $i$ in $\pi$, since all  
  neighbors of $i$ in $\pi$ are larger than $i$.  Furthermore, $i$ is  
  not a neighbor of $i+1$ in $\pi\tau_{n,i-1}$, since $i-1$ is a  
  neighbor of $i+1$ in $\pi\tau_{n,i-1}$, while the second neighbor of  
  $i+1$ in $\pi\tau_{n,i-1}$, if it exists, is larger than $i+1$.  It  
  follows that  
  $$  
  \sigma:=\pi\tau_{n,i-1}\tau_{n,i}\descr\pi.  
  $$  
  Let $\tilde{\pi}:=\pi\tau_{n,i}$, then the same argument as  
  before shows   
  $$  
  \tilde{\sigma}:=\tilde{\pi}\tau_{n,i-1}\tau_{n,i}\descr\tilde{\pi}.  
  $$  
  Furthermore,  
  $$  
  \tilde{\sigma}\tau_{n,i-1} =  
  \tilde{\pi}\tau_{n,i-1}\tau_{n,i}\tau_{n,i-1} =  
  \tilde{\pi}\tau_{n,i}\tau_{n,i-1}\tau_{n,i} =  
  \pi\tau_{n,i-1}\tau_{n,i} = \sigma,  
  $$  
  $i-1$ and $i$ are neighbors in $\sigma$ and each neighbor $\neq  
  i,i-1$ of $i$ and $i-1$ in $\sigma$ is even larger than $i+1$.  
  Inductively,  
  $\pi \descr \sigma \peakr \tilde{\sigma} \descr \tilde{\pi}$.  
\end{proof}

Observe that $\sigma\peakn\pi$ implies that $\Peak(\pi)=\Peak(\sigma)$,  
by Lemma~\ref{dbhl-desc} and \kref{1-2-inv}.  This is thus the immediate  
part of

\begin{lemma} \label{peak-rel}  
  $\Peak(\pi)=\Peak(\sigma)$ if and only if $\pi\peakr \sigma$, for  
  all $\pi,\sigma\in S_n$.  
\end{lemma}  
  
\begin{proof}  
  Let $\pi\in S_n$ and $P:=\Peak(\pi)$. Note that $\Des(\pi)$ is a  
  peak set if and only if $\Des(\pi)=P$. We first show:  
    
  $(*)$ If $\Des(\pi)\neq P$, then there exists a permutation  
  $\tilde{\pi}\in S_n$ such that $\tilde{\pi}\peakr\pi$ and  
  $|\Des(\tilde{\pi})|<|\Des(\pi)|$.  
    
  Let $\Des(\pi)\neq P$, then there exists $k\in\Des(\pi)$ such that  
  $k-1\in\Des(\pi)\cup\{0\}$.  Let $k$ be maximal with this property  
  and put $c:=(k-1)\pi$, $b:=k\pi$ and $i:=(k+1)\pi$. In the case of  
  $k=1$, $c$ should be read as $\infty$. The choice of $k$ implies  
  $c>b>i$.  
    
  The neighbors of $b$ in $\pi$ are $i$ and $c$ (or just $i$ if
  $k=1$).  If $b\neq i+1$, then $b-1$ is not a neighbor of $b$ in
  $\pi$, since $c>b$, hence $\pi\descn\pi\tau_{n,b-1}$.  The neighbors
  of $b-1$ in $\pi\tau_{n,b-1}$ are $c$ and $i$ (or, again, just $i$
  if $k=1$).  Repeating this procedure, gives a permutation
  $\hat{\pi}\descr\pi$ such that $(k-1)\hat{\pi}=c$, $k\hat{\pi}=i+1$
  and $(k+1)\hat{\pi}=i$.  Furthermore,
  $k+1\notin\Des(\pi)=\Des(\hat{\pi})$, which implies $i<(k+2)\pi$ in
  the case of $k+1<n$.  It follows that
  $\tilde{\pi}:=\hat{\pi}\tau_{n,i}\peakr\hat{\pi}$, by
  Proposition~\ref{1-2}, and
  $\Des(\tilde{\pi})=\Des(\pi)\backslash\{k\}$.  This proves $(*)$.
    
  Now let $\sigma\in S_n$ such that $\Peak(\sigma)=\Peak(\pi)=P$.
  Applying $(*)$ a number of times, yields a permutation
  $\tilde{\pi}\in S_n$ such that $\pi\peakr\tilde{\pi}$ and
  $\Des(\tilde{\pi})=\Peak(\tilde{\pi})=P$ and a permutation
  $\tilde{\sigma}\in S_n$ such that $\sigma\peakr\tilde{\sigma}$ and
  $\Des(\tilde{\sigma})=\Peak(\tilde{\sigma})=P$.  By
  Lemma~\ref{dbhl-desc}, 
  $\pi\peakr\tilde{\pi}\descr\tilde{\sigma}\peakr\sigma$.  
\end{proof}

The preceding proposition leads to a  
  
\begin{proof}[Proof of the remaining part of Main Theorem~\ref{comb-char}]  
  Let $\varphi=\sum_{\pi\in S_n} k_\pi \pi\in\!\DD_n$, then
  $k_\pi=k_\sigma$ whenever $\pi\descn\sigma$.  If, additionally,
  $\varphi\tau_{n,1}=\varphi$, it follows that $k_\pi=k_\sigma$
  whenever $\pi\peakn\sigma$, hence $k_\pi=k_\sigma$ for all
  $\pi,\sigma\in S_n$ such that $\pi\peakr\sigma$.  Now
  Lemma~\ref{peak-rel} implies $\varphi\in\PP_n$.
\end{proof}

\breakitdown

The concept of descent and peak relations extends to arbitrary finite  
Coxeter groups as follows.  
  
Let $(W,S)$ be a Coxeter system and assume that $W$ is finite (for  
details, see \cite{solomon76}). $\length(w)$ denotes the length of $w$  
viewed as a word in the elements of $S$, for all $w\in W$.  Let  
$$  
\Des(w) := \setofall{t\in S}{\length(wt)<\length(w)}  
$$  
for all $w\in W$ and denote by $\DD(W)$ the linear subspace of the  
group algebra $KW$ generated by the elements  
$$  
y_K := \sum_{\substack{w\in W\\[2pt] \Des(w)=K}} w,\qquad K\subseteq S.  
$$  
$\tilde{w}$ is called a descent neighbor of $w$ ($\tilde{w}\smile  
w$) if there exists an element $s\in S$ such that $\tilde{w}=sw$ and  
$sw\neq wt$ for all $t\in S$.\footnote{%
  Note that, in the setting of \cite{solomon76} that is used here,  
for $W=S_n$, products of permutations are to be read from right to left.}
Denote by $\sim$ the smallest equivalence relation on $W$ containing  
$\smile$, then  
\begin{equation}  
  \label{desc-rel-cox}  
  \Des(w)=\Des(v)\mbox{ if and only if }w\sim v,  
\end{equation}    
as I learned from discussions with C. Hohlweg and C. Reutenauer.  The
proof is a straight-forward generalization of the proof of
\cite[Proposition~4.2]{blelau93a} and is skipped.
    
As a consequence, there is a short proof of Solomon's result imitating  
the proof of \cite[Theorem~4.1]{blelau93a}.

\begin{theorem*}[Solomon, 1976]  
  $\DD(W)$ is a sub-algebra of $KW$.  
\end{theorem*}  
  
\begin{proof}  
  Let $I,J\subseteq S$, then $y_Iy_J = \sum_{w\in W} \# F(I,J,w)\,w$,
  where
  $$  
  F(I,J,w) := \setofall{(u,v)\in W\times W}{\Des(u)=I,\,\Des(v)=J,\,uv=w},  
  $$  
  hence the claim is equivalent to $ \# F(I,J,w)=\#  
  F(I,J,\tilde{w}) $ for all $w,\tilde{w}\in W$ such that  
  $\Des(w)=\Des(\tilde{w})$, that is $\tilde{w}\sim w$, by  
  \kref{desc-rel-cox}.  
      
  It suffices to consider the case $\tilde{w}\smile w$.  Let $s\in S$  
  such that $\tilde{w}=sw$ and $sw\neq wt$ for all $t\in S$ and define  
  $$  
  \alpha_s:W\times W\lra W\times W,\, 
  (u,v) 
  \lmt  
  \tiptop{%
    (su,v)}{\mbox{if }u^{-1}su\notin S,}{%
    (u,(u^{-1}su)v)}{\mbox{if }u^{-1}su\in S,}  
  $$  
  then $\tilde{u}\tilde{v}=s(uv)$ for all  
  $u,v,\tilde{u},\tilde{v}\in W$ such that  
  $(u,v)\alpha_s=(\tilde{u},\tilde{v})$.  Furthermore,  
  $\tilde{u}^{-1}s\tilde{u}=u^{-1}su$, hence  
  $(\tilde{u},\tilde{v})\alpha_s=(u,v)$.  It follows that $\alpha_s$  
  is bijective and $\alpha_s\alpha_s=\id_{W\times W}$.  
      
  Finally, if $(u,v)\in F(I,J,w)$, then $(\tilde{u},\tilde{v})\in
  F(I,J,\tilde{w})$, which may be derived as follows: If
  $u^{-1}su\notin S$, then $su\neq ut$ for all $t\in S$, which implies
  $\tilde{u}\smile u$ and $\Des(\tilde{u})=\Des(u)=I$.  If
  $\tilde{s}:=u^{-1}su\in S$, then $\tilde{s}v\neq vt$ for all $t\in
  S$, since otherwise $sw=u\tilde{s}v=wt$ for a proper $t\in S$.  It
  follows that $\tilde{v}\smile v$ and $\Des(\tilde{v})=\Des(v)=J$.
      
  Putting together all pieces, $\alpha_s|_{F(I,J,w)}$ is injective,  
  mapping $F(I,J,w)$ into $F(I,J,\tilde{w})$.  Since everything is  
  symmetric in $w$ and $\tilde{w}$, it follows that  
  $$  
  \# F(I,J,w)=\# F(I,J,\tilde{w}),  
  $$  
  and the proof is complete.  
\end{proof}  
    
For any fixed $s^*\in S$, $s^*$-peak relations may be defined by  
$sw\speakn w$ if $sw\smile w$ or $s=s^*$.  
    
The sums of the corresponding $s^*$-peak classes span a right ideal of
the descent algebra $\DD(W)$ --- a \emph{peak-like sub-algebra} of
$\DD(W)$.

%
%
%
%
%
%
%
%
%
%
%
%
%
%
%
%
%
%
%
%
%
%
%
%
%

\section{Series of divided powers} \label{divided}

In the setting of quasi-symmetric functions, Stembridge's peak algebra
occurs as a sub-algebra of $\Qsym$, which, by duality, corresponds to
a homomorphic image of $\DD$.  It turns out that an epimorphism
$(\DD,\conv,\dank)\to(\PP,\conv,\dank)$ of bialgebras is realized by
simultaneous inner right multiplication with $\R^n$ in each
homogeneous component $\DD_n$.  As a consequence, the picture of
duality of $\PP$ and the algebra of peak functions can be completed
(Corollary~\ref{duality-2}).

Let $\alpha_n\in KS_n$ for all $n\ge 0$, then $(\alpha_n)_{n\ge 0}$ is
called a \emph{series of divided powers} if
$$
\alpha_n\dank = \sum_{k=0}^n \alpha_k\tensor\alpha_{n-k},
$$
for all $n\in\N_0$.  Equivalently, the element $\alpha=\sum_{n\ge
  0}\alpha_n$ of the direct product of all $KS_n$, $n\in\N_0$, is
\emph{group-like}, that is
$$
\alpha\dank = \alpha\tensor\alpha.
$$
For instance, $(\Xi^n)_{n\ge 0}$ and $(\R^n)_{n\ge 0}$ are series
of divided powers in $\DD$, by \kref{eq:4} and Main
Theorem~\ref{bialg-P}.

For any $\varphi\in KS$, there exist uniquely determined components
$\varphi_n\in KS_n$ of $\varphi$ such that $\varphi_n=0$ for all but a
finite number of $n$ and $\varphi=\sum_{n\ge 0}\varphi_n$.  Therefore
$\alpha$ yields a linear mapping $KS\to KS$, by (inner) right
multiplication:
$$
\varphi\inn\alpha = \sum_{n\ge 0}\varphi_n\alpha_n\,.
$$
More precisely, there is the following interesting consequence of
the reciprocity law~\ref{reziprozi}:

\begin{lemma} \label{div-pow-epi}
  Let $(\alpha_n)_{n\ge 0}$ be a series of divided powers in $\DD$,
  then right multiplication with $\alpha=\sum_{n\ge 0}\alpha_n$,
  $$
  \iota_\alpha:\DD\lra\DD,\,\varphi\lmt \varphi\inn\alpha,
  $$
  is a graded endomorphism of the bialgebra $(\DD,\conv,\dank)$.
\end{lemma}

\begin{proof}
  Let $\varphi\in\DD_n$, $\psi\in\DD_m$, then
  $\varphi\inn\alpha=\varphi\alpha_n\in\DD_n$, by
  Theorem~\ref{sol-result}, and
  \begin{eqnarray*}
    (\varphi\conv\psi)\inn\alpha
    & = &
    (\varphi\conv\psi)\inn\alpha_{n+m}\\
    & = &
    \Big(
    (\varphi\tensor\psi)\inn\sum_{k=0}^{n+m} \alpha_k\tensor\alpha_{n+m-k}
    \Big)\upnk\\
    & = &
    (\varphi\alpha_n)\conv(\psi\alpha_m)\\[2mm]
    & = &
    (\varphi\inn\alpha)\conv(\psi\inn\alpha),
  \end{eqnarray*}
  by Lemma~\ref{reziprozi}.
  Now choose $\varphi_1^{(k)}\in\DD_k$ and $\varphi_2^{(k)}\in
  \DD_{n-k}$ for all $k\in\haken{n}_0$ such that
  $\varphi\dank=\sum_{k=0}^n\varphi_1^{(k)}\tensor\varphi_2^{(k)}$.
  Since $\dank:\DD\to\DD\tensor\DD$ is a homomorphism for the inner
  product,
  $$
  (\varphi\iota_\alpha)\dank 
  = 
  (\varphi\inn\alpha_n)\dank
  = 
  \varphi\dank\inn\alpha_n\dank
  =
  \sum_{k=0}^n
  \varphi_1^{(k)}\inn\alpha_k\tensor\varphi_2^{(k)}\inn\alpha_{n-k}
  = 
  \varphi\dank\,(\iota_\alpha\tensor\iota_\alpha),
  $$
  which completes the proof, by linearity in $\varphi$ and $\psi$.
\end{proof}

The specialization $\alpha=\R$ gives the main result of this section:

\begin{mtheorem} \label{outer-epi}
  Right multiplication with
  $$
  \R := \sum_{n\ge 0} \R^n
  $$
  yields a graded epimorphism of bialgebras $(\DD,\conv,\dank)\lra
  (\PP,\conv,\dank)$.  In particular,
  $
  \PP_n=\DD_n\R^n
  $
  is a left ideal of $\DD_n$, for all $n\ge 0$.  Furthermore,
  \begin{equation}
    \label{xi-tilde}
    \R^q=\Xi^q\R^n,
  \end{equation}
  for all $n\in\N_0$, $q\zerl n$.
\end{mtheorem}

\begin{proof}
  Right multiplication with $\R$ is an endomorphism of the bialgebra
  $(\DD,\conv,\dank)$, by Lemma~\ref{div-pow-epi}. In particular,
  $$
  \Xi^q\inn\R 
  = 
  (\Xi^{q_1}\conv\cdots\conv\Xi^{q_k})\R^n 
  =
  (\Xi^{q_1}\R^{q_1})\conv\cdots\conv(\Xi^{q_k}\R^{q_k}) 
  =
  \R^q,
  $$
  for all $q=q_1\wldots q_k\zerl n$.  Recalling that
  $\setofall{\Xi^q}{q\in\N^*}$ is a linear basis of $\DD$ and that
  $\setofall{\R^q}{q\in\N^*\mbox{ odd}}$ is a linear basis of $\PP$,
  completes the proof.
\end{proof}

In what follows, $\tilde{\varphi} := \varphi\R$ denotes the image of
$\varphi$ under the epimorphism $\iota_{\R}$, for all $\varphi\in\DD$,
which is compatible with the definition of $\R^q$
given in Section~\ref{peak}, by \kref{xi-tilde}.

Note that both $\DD_n$ and $\PP_n$ are $\DD_n$-left modules (by inner
left multiplication), and --- as a consequence of Main
Theorem~\ref{outer-epi} ---
$$
\DD_n\lra \PP_n, \varphi\lmt \varphi\R^n
$$
is an epimorphism of $\DD_n$-left modules. A more detailed analysis
of the $\DD_n$-left module $\PP_n$ follows in Section~\ref{radical}.

\begin{remark}
  If $(\alpha_n)_{n\ge 0}$ and $(\beta_n)_{n\ge 0}$ are series of
  divided powers in $\DD$, then so is $(\gamma_n)_{n\ge 0}$, defined
  by
  $$
  \gamma_n:=\sum_{k=0}^n \alpha_k\conv\beta_{n-k},
  $$
  for all $n\in\N_0$.  Indeed, 
  $\alpha=\sum_{n\ge 0} \alpha_n$ and $\beta=\sum_{n\ge 0} \beta_n$ 
  are group-like,
  hence also $\gamma=\sum_{n\ge 0} \gamma_n=\alpha\conv\beta$.
  
  The series $(\Delta^{1^{.n}})_{n\ge 0}$ is another series of divided
  powers in $\DD$, as was mentioned in the proof of 
  Main Theorem~\ref{bialg-P}
  already.  Putting $\Delta:=\sum_{n\ge 0} \Delta^{1^{.n}}$,
  $$
  \R=\Delta\conv \Xi.
  $$
  It would be interesting to apply Proposition~\ref{div-pow-epi} to
  other monomials $\alpha$ in $\Delta$ and $\Xi$ (instead of $\R$) and
  to analyze the corresponding homomorphic images $\DD\iota_\alpha$ of
  $\DD$, each of which is at the same time a sub-bialgebra of
  $(\DD,\conv,\dank)$ and a left ideal of $(\DD,\inn)$.
\end{remark}

Repeated right multiplication with $\R$ still gives $\PP$, as is derived now.
First two helpful formulae:

\begin{proposition} \label{multi-tilde-xi}
  $
  \Delta^{1^{.n}}\R^n=\R^n
  $
  and 
  $
  \R^n\R^n = 2\R^n + \sum_{k=1}^{n-1} \R^{k.(n-k)}
  $,
  for all $n\in\N_0$.
\end{proposition}

\begin{proof}
  Let $n\in\N_0$ and $k\in\haken{n-1}_0$, then
  $\Delta^{1^{.n}}\Delta^{1^{.k}.(n-k)}=\Delta^{1^{.(n-k-1)}.(k+1)}$, since
  $\Delta^{1^{.n}}$ is the order reversing involution in $S_n$.
  Equation \kref{re-xi-tilde} implies
  \begin{eqnarray*}
    \Delta^{1^{.n}}\R^n
    & = &
    2\sum_{k=0}^{n-1} \Delta^{1^{.n}}\Delta^{1^{.k}.(n-k)}
      = 
    2\sum_{k=0}^{n-1} \Delta^{1^{.(n-k-1)}.(k+1)}
      = 
    \R^n
  \end{eqnarray*}
  as asserted. In particular,
  $$
  \R^n\R^n
  =
  \sum_{k=0}^n (\Delta^{1^{.k}}\R^k)\conv(\Delta^{n-k}\R^{n-k})
  =
  2\R^n + \sum_{k=1}^{n-1} \R^{k.(n-k)},
  $$
  by \kref{re-xi-tilde} and Main Theorem~\ref{outer-epi}.
\end{proof}

\begin{corollary}
  $\PP_n\R^n = \PP_n\,$, for all $n\in\N_0$.  In particular, right
  multiplication with $\R$ yields an automorphism of the bialgebra
  $(\PP,\conv,\dank)$.
\end{corollary}

\begin{proof}
  Let $q=q_1\wldots q_k\zerl n$, then 
  $$
  \R^q\R^n
  =
  (\R^{q_1}\R^{q_1})\conv\cdots\conv(\R^{q_k}\R^{q_k})
  \in
  2^k\R^q + \erz{\setofall{\R^r}{r\zerl q,\,\length(r)>\length(q)}}_K\,,
  $$
  by Main Theorem~\ref{outer-epi} and
  Proposition~\ref{multi-tilde-xi}.  The elements $\R^q\R^n$,
  $q\zerlodd n$, are thus linearly independent, by
  Proposition~\ref{transition-tilde-xi}, and both claims follow.
\end{proof}

\breakitdown

Aiming at the dual version of Main Theorem~\ref{duality-1}, a proper
description of the kernel of $\iota_{\R}$ is needed. This description
is based on the following simple, but lengthy computation, which will
also be useful in Section~\ref{PLI}.  Recall that $D\sctriangle(D+1)$
denotes the symmetric difference of the sets $D$ and $D+1$.

\begin{proposition} \label{tilde-delta}
  For all $D\subseteq \haken{n-1}$,
  $$
  \tilde{\Delta}^D = \Delta^D\R^n = \sum_P 2^{|P|+1} \Pi^P,
  $$
  summed over all peak sets $P$ in $\haken{n}$ such that $P\subseteq
  D\sctriangle(D+1)$.
\end{proposition}

\begin{proof}
  By \kref{eq:1} and Proposition~\ref{transition-tilde-xi},
  \begin{eqnarray*}
    \tilde{\Delta}^D
    & = &
    \sum_{E\subseteq D} (-1)^{|D|-|E|}\; \Xi^E\R^n\\
    & = &
    \sum_{E\subseteq\haken{n-1}}
    \sum_{E\subseteq D} (-1)^{|D|-|E|}\; 2^{|E|+1}
    \sum_{P\subseteq E\cup (E+1)}\;
    \Pi^P\\
    & = &
    \sum_P
    \sum_{\substack{E\subseteq D\\[2pt] P\subseteq E\cup (E+1)}} 
    (-1)^{|D|-|E|}\; 2^{|E|+1}
    \Pi^P.
  \end{eqnarray*}
  Let $P$ be a peak set in $\haken{n}$ and put
  $$
  k_{D,P} := \sum_{\substack{E\subseteq D\\[2pt] P\subseteq E\cup (E+1)}}
  (-1)^{|D|-|E|}\; 2^{|E|+1},
  $$
  then $k_{D,P}\neq 0$ implies that $P\subseteq D\cup (D+1)$.  If
  $P=\emptyset$, then
  $$
  k_{D,P}
  = 
  (-1)^{|D|} \sum_{E\subseteq D} (-1)^{|E|}\; 2^{|E|+1} 
  = 
  2\sum_{j=0}^{|D|} \binom{|D|}{j}\,(-1)^{|D|-j}\; 2^j 
  = 
  2
  $$
  as asserted.

  Let $P\neq\emptyset$ and choose $i\in P$, then $i-1\in E$ or $i\in
  E$, for each $E\subseteq\haken{n-1}$ such that $P\subseteq E\cup
  (E+1)$. Set $\tilde{P}:=P\backslash\{i\}$ and
  $\tilde{D}:=D\backslash\{i-1,i\}$, then
  \begin{eqnarray*}
    k_{D,P}
    & = &
    \sum_{\substack{E\subseteq D\\[2pt] P\subseteq E\cup (E+1)}} 
    (-1)^{|D|-|E|} 2^{|E|+1}\\
    & = &
    (-1)^{|D|}
    \Bigg[
    \sum_{\substack{i\notin E\subseteq D\\[2pt] P\subseteq E\cup (E+1)}}
    (-1)^{|E|}\; 2^{|E|+1}
    +
    \sum_{\substack{i-1\notin E\subseteq D\\[2pt] P\subseteq E\cup (E+1)}}
    (-1)^{|E|}\; 2^{|E|+1}\\
    & + &
    \sum_{%
     \substack{\{i-1,i\}\subseteq E\subseteq D\\[2pt] P\subseteq E\cup (E+1)}} 
    (-1)^{|E|}\; 2^{|E|+1}
    \Bigg]\\
  \end{eqnarray*}
  But, if $i-1,i\in D$, it follows that
  $$
  k_{D,P} 
  = 
  (-1)^{|D|} 
  \sum_{\substack{E\subseteq \tilde{D}\\[2pt] \tilde{P}\subseteq E\cup (E+1)}}
  (-1)^{|E|}\; 2^{|E|+1}
  \Big[-2-2+4\Big] = 0,
  $$
  while, if $|D\cap\{i-1,i\}|=1$,
  $$
  k_{D,P} 
  = 
  (-1)^{|D|} 
  \sum_{\substack{E\subseteq \tilde{D}\\[2pt] \tilde{P}\subseteq E\cup (E+1)}}
  (-1)^{|E|}\; 2^{|E|+1}
  \Big[-2\Big] 
  = 
  2k_{\tilde{D},\tilde{P}}.
  $$
  Now, inductively, $k_{D,P}\neq 0$ if and only if either $i\in D$ or
  $i-1\in D$, for all $i\in P$, and in this case
  $k_{D,P}=2^{|P|}k_{D,\emptyset}=2^{|P|+1}$.
  The claim follows by observing that the condition that either
  $i\in D$ or $i-1\in D$, for all $i\in P$, is equivalent to
  $P\subseteq D\sctriangle (D+1)$.
\end{proof}

Recall the definition of the pairing $[\,\cdot\,,\,\cdot\,]$ given in
Main Theorem~\ref{duality-1} to state:

\begin{corollary} \label{duality-2}
 $\ker\,\iota_{\R}=\PPP^\perp$.

 In particular, $\PPP$ is a sub-coalgebra of $(\Qsym,\gamma)$, dual to
 the algebra $(\PP,\conv)$, and
 $$
 [\varphi\conv\psi,f]
 =
 [\varphi\tensor\psi,\gamma(f)]
 $$
 for all $\varphi,\psi\in\PP$, $f\in\PPP$.
\end{corollary}

\begin{proof}
  Let $\varphi=\sum_{E\subseteq\haken{n-1}} a_E\Delta^E\in\DD_n$, then
  $$
  0
  =
  <\varphi,K_P>
  =
  \sum_{E\subseteq\haken{n-1}} \sum_{P\subseteq D\scsctriangle(D+1)} 
  a_E <\Delta^E,F_D>
  =
  \sum_{P\subseteq D\scsctriangle(D+1)} a_D
  $$
  for all peak sets $P$ in $\haken{n}$ is equivalent to
  $$
  \tilde{\varphi}
  =
  \sum_E a_E\tilde{\Delta}^E
  =
  \sum_E a_E\sum_{P\subseteq D\scsctriangle(D+1)} 
  \Pi^P
  =
  \sum_P
  \Big(\sum_{P\subseteq D\scsctriangle(D+1)} a_D\Big)\Pi^P
  =
  0.
  $$
  This proves the first claim.
  All what remains is now obtained by exchanging the roles of $\PP$ and $\PPP$
  and replacing $\theta$ by $\iota_{\R}$ in the proof of
  Main Theorem~\ref{duality-1}.
\end{proof}

%
%
%
%
%
%
%
%
%
%
%
%
%
%
%
%
%
%
%
%
%
%
%
%
%

\section{Solomon's epimorphism}
\label{c}

The descent algebra $\DD_n$ is linked to the character theory of the
symmetric group $S_n$ by means of an epimorphism of algebras
$$
\DD_n\lra \Cl_K(S_n)
$$
introduced by Solomon (\cite{solomon76}).  Here $\Cl_K(S_n)$
denotes the ring of $K$-valued class functions of $S_n$.  In this
section the image of $\PP_n$ under Solomon's epimorphism is
determined.  This result serves as a helpful tool for the
investigations in inner products in $\PP$ in the sections that follow.

Some notations are needed.
Let $C_q$ be the conjugacy class consisting of all permutations $\pi$ whose
cycle partition is a rearrangement of $q$, for all $q\in\N^*$.  If $a_i$
denotes the multiplicity of the letter $i$ in $q$, for all $i\in\N$ then
$$
q?
:=
\prod_{i\ge 1} a_i!\;i^{a_i}
$$
is the order of the centralizer of any $\pi\in C_q$ in $S_n$.  Let
$\cha_q$ be the characteristic function of $C_q$ in $\Cl_K(S_n)$,
mapping each $\pi\in S_n$ to $1$ if $\pi\in C_q$, and to $0$
otherwise; and set $\ch_q:=q?\cha_q$.  To save trouble, let it be
mentioned that $C_q=C_r$, $q?=r?$, $\cha_q=\cha_r$ and $\ch_q=\ch_r$
whenever $q$ is a rearrangement of $r$ ($q\ass r$).

If $p=p_1\wldots p_l\zerl n$ and $p_1\ge \cdots \ge p_l$, then $p$ is
called a \emph{partition} of $n$, and we write $p\p n$. If, additionally,
$p$ is odd, write $p\podd  n$. Set
$$
\Part(n)
:= \setofall{p}{p\p n},
$$
then the elements $\cha_p$, $p\in\Part(n)$, constitute a linear
basis of $\Cl_K(S_n)$.
The direct sum
$$
\CC:=\bigoplus_{n\in\N_0} \Cl_K(S_n)
$$
is equipped with an \emph{inner product} $\cdot$ inherited from the ordinary
products of $K$-valued functions in $\Cl_K(S_n)$ --- sometimes
referred to as the \emph{Kronecker product} --- by orthogonal extension:
$$
\chi\cdot\psi:=\tiptop{\chi\psi}{\mbox{if $m=n$}}{0}{\mbox{otherwise}}
$$
for all $n,m\in\N$, $\chi\in\Cl_K(S_n)$, $\psi\in\Cl_K(S_m)$.

The \emph{outer product} $\bullet$ on $\CC$ corresponds to the
ordinary multiplication of symmetric functions, via Frobenius'
characteristic mapping, and may be defined by
\begin{equation}
  \label{chqchr}
  \ch_r\bullet\ch_q=\ch_{r.q},
\end{equation}
for all $q,r\in\N^*$.  For all $n\in\N_0$, let $\xi^n$ denote the
trivial character of $S_n$.  By Theorem~\ref{bialg-D}, there is a
unique homomorphism of algebras $c:(\DD,\conv)\to (\CC,\bullet)$ such
that
$$
c(\Xi^n) = \xi^n
$$
for all $n$.  Another crucial result due to Solomon
(\cite{solomon76}) is the following.

\begin{theorem}\label{sol-rad}
  $c:\DD\to \CC$ is an epimorphism of algebras with respect to the
  \emph{inner} products, and
  $$
  \ker\,c|_{\DD_n} = \Big\langle\,
  \Setofall{%
    \Xi^q-\Xi^{\dot{q}}}{q,\dot{q}\zerl n,\, q\ass \dot{q} }
  \,\Big\rangle_K
  $$
  for all $n\in\N_0$.
\end{theorem}

Let
$$
\CC^{odd} 
:= 
\Big\langle\,\Setofall{\cha_q}{q\in\N^*\mbox{ odd}}\,\Big\rangle_K\,,
$$
then the dimension of the $n$-th homogeneous component of $\CC^{odd}$ is the
number of odd partitions of $n$.

\begin{mtheorem} \label{p-c-im}
  $c|_\PP:\PP\to \CC^{odd}$ is an epimorphism of algebras for both the
  inner and the outer products, and
  $$
  \ker\,c|_{\PP_n} = \Big\langle\,
  \Setofall{%
    \R^q-\R^{\dot{q}}}{q,\dot{q}\zerlodd  n,\, q\ass \dot{q} }
  \,\Big\rangle_K
  $$
  for all $n\in\N_0$. Furthermore,
  $
  c(\R^n)
  = 
  \sum_{p\podd\,  n} 2^{\length(p)}\,\cha_p
  $.
\end{mtheorem}

\begin{proof}
  The image of $\Delta^{1^{.k}}$ under $c$ is the the sign character
  $\eps_k$ of $S_k$, for all $k\in\N_0$, hence \kref{re-xi-tilde}
  implies that, for odd $n$,
  \begin{eqnarray*}
    c(\R^n)
    & = &
    \eps_n
    +
    \sum_{k=1}^{n-1} (\eps_k\bullet\xi_{n-k})
    +
    \xi^n
    = 
    \sum_{p\podd  n} 2^{\length(p)}\,\cha_p
  \end{eqnarray*}
  (see, for instance, \cite[III.8, Example~6(a)]{macdonald95}). As a 
  consequence, $c(\PP)=\CC^{odd}$, by Theorem~\ref{sol-rad} and Main
  Theorem~\ref{outer-epi}.
  Furthermore, $c(\R^q-\R^{\dot{q}})=0$ whenever $q\ass \dot{q}$, since
  $(\CC,\bullet)$ is commutative.  Comparing dimensions, completes the proof.
\end{proof}

Note that Main Theorem~\ref{p-c-im} is the dual way of stating 
\cite[Theorem~3.8]{stembridge97}.

%
%
%
%
%
%
%
%
%
%
%
%
%
%
%
%
%
%
%
%
%
%
%
%
%

\section{Peak Lie idempotents} \label{PLI}

Each of the classical Lie idempotents --- the Dynkin operator, the
canonical Lie idempotent, and the Klyachko idempotent --- is contained
in $\DD$.  The study of these Lie idempotents played a major role
in the investigation of combinatorial as well as algebraic
properties of the Solomon algebra (see, for instance,
\cite{bauer01,blelau96,garsia-reutenauer89}).  Apart from that, the
theory of Lie idempotents is interesting for its own sake
and allows numerous combinatorial applications
(\cite{bbg90,garsia89,LST96,patrasreut99,schocker-derange}). 

In this section, the basics concerning Lie idempotents contained in
$\PP$ --- \emph{peak Lie idempotents} --- are discussed.  As a first
consequence, more detailed information on the structure of
$(\PP,\conv,\dank)$ and its primitive Lie algebra is obtained (Main
Theorem~\ref{high-lie-id-basis}).  Furthermore, peak variants of the
classical Lie idempotents are displayed
(Proposition~\ref{peak-dynkin}, Proposition~\ref{peak-canon},
\kref{peak-klyachko}).  In the case of the peak Dynkin operator, two
combinatorial applications are given exemplarily
(Corollary~\ref{comb-cor-1}, Corollary~\ref{comb-cor-2}), while a peak
variant of the algebra of noncommutative cyclic characters
(\cite{LST96}) arises from the Klyachko peak Lie idempotent
(Corollary~\ref{peak-noncyclic}).
 
To begin with, recall the definition of a Lie idempotent. Let $X$ be
an infinite alphabet and denote by $A(X)$ the free associative algebra
over $X$.  Let $A_n(X)$ be the $n$-th homogeneous component of $A(X)$,
for all $n\in\N_0$, then $S_n$ acts on $A_n(X)$ from the left, via
Polya action:
$$
\pi(x_1\cdots x_n) := x_{1\pi}\cdots x_{n\pi}
$$
for all $\pi\in S_n$, $x_1,\ldots,x_n\in X$.

The usual Lie product $a\circ b=ab-ba$ defines the structure of a Lie
algebra on $A(X)$, and the Lie sub-algebra $L(X)$ generated by $X$ is
free over $X$ by a classical result of Witt (\cite{witt37}).  Let
$L_n(X):=L(X)\cap A_n(X)$ denote the $n$-th homogeneous component of
$L(X)$, for all $n\in\N_0$.

An element $\alpha\in KS_n$ is called an idempotent, if
$\alpha^2=\alpha$.  Due to Dynkin (\cite{dynkin47}), Specht (\cite{specht48})
and Wever (\cite{wever49}),
$$
\omega_n := \sum_{k=0}^{n-1}(-1)^k\Delta^{1^{.k}.(n-k)}
$$
is (up to the factor $\frac{1}{n}$) an idempotent in $\DD_n$,
sometimes called the \emph{Dynkin operator}.  Via Polya action,
\begin{equation}
  \label{left-normed}
  \omega_n\,(x_1\cdots x_n)
  = (\cdots((x_1\circ x_2)\circ x_3)\circ\quad\cdots\quad)\circ x_n\,,  
\end{equation}
for all $x_1,\ldots,x_n\in X$, hence $\frac{1}{n}\omega_n$ projects
$A_n(X)$ onto $L_n(X)$.  Accordingly, any idempotent $\alpha\in KS_n$ is
called a \emph{Lie idempotent} if $\alpha KS_n=\omega_n KS_n$.

Note that, if $\alpha,\beta\in KS_n$ and $\alpha$ is a Lie idempotent,
then $\beta$ is a Lie idempotent if and only if $\alpha\beta=\alpha$
and $\beta\alpha=\alpha$.  Furthermore, for any Lie idempotent
$\alpha\in KS_n$ and any $q\zerl n$,
\begin{equation}
  \label{xiq-lie-id}
  \Xi^q\alpha
  =
  \textstyle
  \frac{1}{n}\Xi^q\omega_n\alpha
  =
  0
  \mbox{ whenever }
  \length(q)>1
\end{equation}
(\cite[Theorem~1.2 and its proof]{bbg90}).
Lie idempotents in the peak algebra are constructed easily:

\begin{proposition} \label{ex-peak-lie}
  Let $n\in\N$ and $\alpha\in\DD_n$ be a Lie idempotent, then
  $$
  \R^n\alpha = \tiptop{2\alpha}{\mbox{if $n$ is
      odd,}}{0}{\mbox{otherwise.}}
  $$
  In particular, if $n$ is odd, then
  $
  \frac{1}{2}\tilde{\alpha} = \frac{1}{2}\alpha\R^n
  $
  is a Lie idempotent in $\PP_n$.
\end{proposition}

\begin{proof}
  By Proposition~\ref{transition-tilde-xi} and \kref{xiq-lie-id},
  \begin{eqnarray*}
    \R^n\alpha
    & = &
    2
    \sum_{\substack{q\zerl n\\[2pt] q^\dagger\;\mbox{\tiny odd}}}
    (-1)^{n-\length(q)}\;
    \Xi^q\alpha\\
    & = &
    \tiptop{2(-1)^{n-1}\;\Xi^n\alpha}{\mbox{if $n$ is odd,}}
           {                     0}{\mbox{otherwise,}}\\
    & = &
    \tiptop{2\alpha}{\mbox{if $n$ is odd,}}
           {      0}{\mbox{otherwise.}}
  \end{eqnarray*}
  In particular,
  $\alpha\tilde{\alpha}=\alpha(\alpha\R^n)=\alpha\R^n=\tilde{\alpha}$
  and $\tilde{\alpha}\alpha=\alpha(\R^n\alpha)=2\alpha^2=2\alpha$, if
  $n$ is odd.
\end{proof}

Any Lie idempotent $\alpha\in\PP$ is called a \emph{peak Lie
  idempotent}.
Recalling that
\begin{equation}
  \label{c-von-lieid}
  c(\alpha)=\cha_n,
\end{equation}
for each Lie idempotent $\alpha\in\DD_n$
(\cite[Proposition~1.8]{blelau02}) gives:

\begin{corollary} \label{when-peak-lie}
  Let $n\in\N$, then $\PP_n$ contains a Lie idempotent if and only if
  $n$ is odd.

  Furthermore, if $n$ is odd and $\alpha\in\DD_n$ is a Lie idempotent,
  then the linear span of all peak Lie idempotents in $KS_n$ is
  $\alpha\PP_n$.  
\end{corollary}

\begin{proof}
  If $n$ is odd, then $\PP_n$ contains the Lie idempotent
  $\frac{1}{2n}\tilde{\omega}_n$, by Proposition~\ref{ex-peak-lie}.  If,
  conversely, $\alpha\in\PP_n$ is a peak Lie idempotent, then
  $\cha_n=c(\alpha)\in c(\PP_n)$, by \kref{c-von-lieid}, hence $n$ is
  odd, by Main Theorem~\ref{p-c-im}.
  
  The linear span of all Lie idempotents in $\DD_n$ is $\omega_n\DD_n$
  (\cite[Proposition~4.5]{blelau93a}). Let $n\in\N$ be odd and choose
  an arbitrary Lie idempotent $\alpha\in\DD_n$, then the defining
  identities $\omega_n\alpha=n\alpha$ and $\alpha\omega_n=\omega_n$
  imply $\omega_n\DD_n=\alpha\DD_n$.  It follows that the linear span
  of all peak Lie idempotents in $KS_n$ is
  $\alpha\DD_n\cap\PP_n=\alpha\PP_n$ as asserted.
\end{proof}

Let $\gamma:\N\lra\DD$ such that $\frac{1}{n}\gamma_n$ is a Lie
idempotent, for all $n\in\N$, then $\gamma$ is called a \emph{Lie
  series} in $\DD$.  Set $ \gamma_q =
\gamma_{q_1}\conv\cdots\conv\gamma_{q_k} $ for all $q=q_1\wldots
q_k\in\N^*$.  If $\gamma$ and $\delta$ are Lie series in $\DD$ and
$q,r\in\N^*$, then
\begin{equation}
  \label{eq:pr}
  \delta_r\gamma_q=r?\gamma_r\mbox{ whenever }r\ass q
\end{equation}
(\cite[Lemma 3.1(ii)]{patrasreut99}). In particular,
$\frac{1}{q?}\gamma_q$ is an idempotent, for all $q\in\N^*$.

In the case $X=\N$, $A(X)=A(\N)$ is the semi-group algebra of $\N^*$ over $K$,
and the mapping $\gamma:\N\to\DD$ extends to an isomorphism of algebras
\begin{equation}
  \label{gamma-iso}
  \gamma:A(\N)\lra(\DD,\conv),\,\sum_qk_qq\lmt \sum_qk_q\gamma_q\,.
\end{equation}
(\cite[Theorem~5.1]{patrasreut99}).  Denote by $\N_{odd}$ the set of
all odd $n\in\N$, then $A(\N_{odd})$ may be viewed as a sub-algebra of
$A(\N)$.  Recall that, if $(B,\cdot,\delta)$ is a bialgebra with
identity $1$, then $b\in B$ is \emph{primitive} if $b\delta=b\tensor
1+1\tensor b$.  The set $\Prim(B)$ of primitive elements in $B$ forms
a Lie sub-algebra of the Lie algebra associated to $(B,\cdot)$.  Each
element of $\omega_nKS_n$ is primitive in $(KS,\conv,\dank)$
(see \cite{patrasreut01}).

\begin{mtheorem} \label{high-lie-id-basis}
  Let $(\gamma_n)_{n\in\N}$ be a Lie series in $\DD$ such that
  $\gamma_n\in\PP_n$ for all odd $n$, then
  $$
  \setofall{\frac{1}{q?}\gamma_q}{q\in\N^*\mbox{ odd}}
  $$
  is a linear basis of $\PP$ consisting of idempotents, and $\gamma$ yields
  an isomorphism of algebras
  $$
  A(\N_{odd})\lra (\PP,\conv),
  $$
  by restriction. Furthermore, the primitive Lie algebra of
  $(\PP,\conv,\dank)$ is 
  $$
  \Prim(\PP)
  =
  \bigoplus_{n\ge 0} \gamma_n\PP_n\,.
  $$
  This Lie algebra is free, freely generated by the elements
  $\gamma_n$, $n\in\N_{odd}$, and $\PP$ is its enveloping algebra.  The
  dimension of the $n$-th homogeneous component $ \gamma_n\PP_n$
  of $\Prim(\PP)$ is equal to the number of Lyndon words $q\in\N^*$
  such that $q\zerlodd n$.
\end{mtheorem}
(Concerning the definition and basic properties of Lyndon words, the
reader is referred to \cite[Section~3]{garsia89}.)

Note that, for odd $n$, the last part of the preceding result gives
the dimension of the linear span of all peak Lie idempotents in
$KS_n$, by Corollary~\ref{when-peak-lie}.  As another consequence,
$(\PP,\conv,\dank)$ is a concatenation Hopf algebra in the
non-commuting variables $\gamma_n$, $n\in\N_{odd}$ (see
\cite[p.~969]{malvenuto-reutenauer95}). This result, however, may also
be obtained directly from Main Theorem~\ref{bialg-P} (ibid., part 3.
of the proof of Theorem~2.1).

\begin{proof}
  We point out that all the assertions hold true if $\PP$ is replaced
  by $\DD$ and the word \emph{odd} is removed.  This has now become
  folklore, although --- as it seems --- some parts are contained in
  the literature only implicitly.  The properties of $\PP$ stated
  above are immediate consequences of the corresponding properties of
  $\DD$.  For the reader's convenience, an explicit proof follows.

  First, by \kref{eq:pr} and \kref{gamma-iso},
  $\setofall{\frac{1}{q?}\gamma_q}{q\in\N^*}$
  is a linear basis of $\DD$ consisting of idempotents. Comparing
  dimensions, implies the first claim. In particular, $\gamma$
  yields an isomorphism of algebras $A(\N_{odd})\to\PP$ when
  restricted to $A(\N_{odd})$ as asserted.
  
  If $\varphi\in\DD_n$ is primitive, then $\Xi^q\varphi=0$ for all
  $q\zerl n$ such that $\length(q)>1$, by Lemma~\ref{reziprozi}, hence
  $\varphi=\frac{1}{n}\omega_n\varphi\in\omega_n\DD_n$, by
  \cite[Proposition~1.2]{blelau96}.  Conversely, any element of
  $\omega_n\DD_n\subseteq \omega_nKS_n$ is primitive, as was mentioned
  above, hence
  $$
  \Prim(\DD)
  =
  \bigoplus_{n\ge 0} \omega_n\DD_n\,.
  $$
  In particular, the set of primitive elements in
  $\PP_n$ is $\omega_n\DD_n\cap\PP_n=\omega_n\PP_n=\gamma_n\PP_n$.
  
  As is readily seen from \kref{left-normed}, the elements $\omega_kq$
  ($q\in\N^*$ odd, $k=\length(q)$) linearly span $L(\N_{odd})$.  The
  image of $L(\N_{odd})$ under $\gamma$ is thus
  $$
  L(\N_{odd})\gamma
  =
  \Erz{\Setofall{\gamma_{\omega_kq}}{q\in\N^*\mbox{ odd},\quad k=\length(q)}}_K
  \,.
  $$
  But \kref{eq:pr} and the multiplication rule
  \cite[Theorem~1.5(b)]{blelau96} imply
  $$
  \omega_n\gamma_q
  =
  \textstyle
  \frac{1}{q?} \omega_n\omega_q\gamma_q
  =
  q_1\frac{1}{q?} \omega_{\omega_kq}\gamma_q
  =
  q_1\gamma_{\omega_kq}\,,
  $$
  thus
  $$
  \Prim(\PP)
  =
  \bigoplus_{n\ge 0} \omega_n\PP_n
  =
  L(\N_{odd})\gamma.
  $$
  All the remaining assertions now follow from the fact the $L(\N_{odd})$ is
  free over $\N_{odd}$.
\end{proof}

For the remainder of this section, the peak variants of the classical
Lie idempotents are studied.

\begin{center}
  \textbf{The peak Dynkin operator}
\end{center}
Let $n$ be odd, then
$
\tilde{\omega}_n = \omega_n\R^n
$
is (up to the factor $\frac{1}{2n}$) a peak Lie idempotent, by
Proposition~\ref{ex-peak-lie}.  A short calculation yields:

\begin{proposition} \label{peak-dynkin}
  $
  \frac{1}{2n}\tilde{\omega}_n 
  = 
  \frac{1}{n} 
  \Big(
  \Pi^{\emptyset} + 2\sum_{k=2}^{n-1}(-1)^{k-1}\Pi^{\{k\}} 
  \Big)
  $,
  for all odd $n$.
\end{proposition}

\begin{proof}
  For all $k\in\haken{n-1}_0$,
  $$
  \tilde{\Delta}^{\haken{k}}
  =
  \sum_{P\subseteq\{1,k+1\}} 2^{|P|+1} \Pi^P
  =
  \tiptop{2\Pi^\emptyset               }{\mbox{if $k\in\{0,n-1\}$,}}
         {2\Pi^\emptyset+4\Pi^{\{k+1\}}}{\mbox{otherwise,}}
  $$
  by Proposition~\ref{tilde-delta}, hence for odd $n$,
  \begin{eqnarray*}
    \tilde{\omega}_n 
    & = & 
    \sum_{k=0}^{n-1} (-1)^k   \tilde{\Delta}^{1^{.k}.(n-k)}\\
    & = & 
    2\sum_{k=0}^{n-1} (-1)^k   \Pi^\emptyset
    +
    4\sum_{k=1}^{n-2} (-1)^k   \Pi^{\{k+1\}}\\
    & = & 
    2\Pi^\emptyset
    +
    4\sum_{k=2}^{n-1} (-1)^{k-1} \Pi^{\{k\}}
  \end{eqnarray*}
  as asserted.
\end{proof}

Two combinatorial applications of the preceding result follow.

Fix $n\in\N$ and denote by $ \P_0 $ the set of all permutations
$\pi\in S_n$ without peak.  Furthermore, let $\P_{1,e}$ (respectively,
$\P_{1,o}$) be the set of all $\pi\in S_n$ with a single peak $k$,
which is even (respectively, odd).

Any partition $p=p_1\wldots p_k$ is called a \emph{block}, if
$p_1=p_2=\cdots=p_k$. Denote by $\mu$ the number theoretic M\"obius
function and define
$$
\sign(p) := \tiptop{\mu(p_1)}{\mbox{if $p$ is a
    block,}}{0}{\mbox{otherwise.}}
$$
 
\begin{corollary} \label{comb-cor-1}
  Let $n$ be odd and $p\p n$, then
  $$
  2|C_p\cap \P_{1,e}|+\sign(p) = 2|C_p\cap \P_{1,o}|+|C_p\cap
  \P_0|.
  $$
  In particular, the number of permutations $\pi$ of cycle type $p$
  in $S_n$ without peak is even if and only if $p$ is not a block or
  $p_1$ is divisible by a square.
\end{corollary}

\begin{proof}
  If $\alpha=\sum_{\pi\in S_n} k_\pi \pi\in KS_n$ such that
  $\frac{1}{n}\alpha$ is a Lie idempotent, then
  $$
  \sum_{\pi\in C_p} k_\pi = \sign(p)
  $$
  (\cite[Proposition~5.1]{garsia89}).  Applied to
  $\alpha=\frac{1}{2}\tilde{\omega}_n$, this proves the claim.
\end{proof}

To push things deeper into that direction, let $\T^n$ be the set of
all standard tableaux with $n$ entries and denote the set of all $t\in
\T^n$ with $l$ odd columns by $\T^{n,l}$, for all $l\in\N_0$. For
instance,
$$
\T^{3,1}
=
\Bigg\{\;
\raisebox{-5mm}{$\unitlength4mm
\linethickness{0.4pt}
\begin{picture}(1.00,3.00)
\put(  0.00, 0.00){\Kasten{1}}

\put(  0.00, 1.00){\Kasten{2}}
\put(  0.00, 2.00){\Kasten{3}}
\end{picture}$}\quad,\;
\raisebox{-5mm}{$\unitlength4mm
\begin{picture}(2.00,2.00)
\put(  0.00, 0.00){\Kasten{1}}

\put(  1.00, 0.00){\Kasten{2}}
\put(  0.00, 1.00){\Kasten{3}}
\end{picture}$}\quad,\;
\raisebox{-5mm}{$\unitlength4mm
\begin{picture}(2.00,2.00)
\put(  0.00, 0.00){\Kasten{1}}

\put(  1.00, 0.00){\Kasten{3}}
\put(  0.00, 1.00){\Kasten{2}}
\end{picture}$}\quad
\Bigg\}
\quad\mbox{ and }\quad
\T^{3,3}
=
\Bigg\{\;
\raisebox{-1mm}{$\unitlength4mm
\linethickness{0.4pt}
\begin{picture}(3.00,1.00)
\put(  0.00, 0.00){\Kasten{1}}

\put(  1.00, 0.00){\Kasten{2}}
\put(  2.00, 0.00){\Kasten{3}}
\end{picture}$}\quad
\Bigg\},
$$
while $\T^{3,2}=\emptyset$.
Let $k\in\{2,\ldots,n-1\}$ and $t\in \T^n$, then $k$ is called a
\emph{peak} of $t$ if $k+1$ stands strictly above $k$ in $t$, while
$k-1$ does not.  As in the case of permutations, let $\T^{n,l}_{1,e}$
(respectively, $\T^{n,l}_{1,o}$) be the set of all $t\in\T^{n,l}$ with
exactly one peak, which is even (respectively, odd).
For example,
$$
\T^{3,1}_{1,e}
=
\Bigg\{
\raisebox{-3mm}{$\unitlength4mm
\begin{picture}(2.00,2.00)
\put(  0.00, 0.00){\Kasten{1}}

\put(  1.00, 0.00){\Kasten{2}}
\put(  0.00, 1.00){\Kasten{3}}
\end{picture}$}\quad
\Bigg\}
\quad\mbox{ and }\quad
\T^{3,1}_{1,o}
=
\emptyset.
$$
Finally, let
$
\T^n_{1,e}
:=
\bigcup_{l\ge 0} \T^{n,l}_{1,e}
$ and
$
\T^n_{1,o}
:=
\bigcup_{l\ge 0} \T^{n,l}_{1,o}
$.

\begin{corollary} \label{comb-cor-2}
  Let $n\in\N$ be odd and $l\in\haken{n}_0$, then
  $$
  |\T^{n,l}_{1,e}|=|\T^{n,l}_{1,o}|+1,
  $$
  if $l$ is odd and $l<n$, while otherwise,
  $\T^{n,l}_{1,e}=\emptyset=\T^{n,l}_{1,o}$.  In particular,
  $$
  |\T^n_{1,e}|=|\T^n_{1,o}|+(n-1)/2.
  $$
\end{corollary}

\begin{proof}
  To start with, observe that $\T^{n,l}\neq\emptyset$ implies
  that $l$ is odd, since $n$ is odd.
  Furthermore, $\T^{n,n}$ consists of the single tableau 
  $
  t
  =
  \raisebox{-1mm}{$\unitlength4mm
    \begin{picture}(6.50,1.00)
      \put(  0.00, 0.00){\Kasten{1}}
      \put(  1.00, 0.00){\Kasten{2}}
      \put(  5.00, 0.00){\Kasten{$n$}}
      \put(  3.25, 0.20){\mbox{$\cdots$}}
      \put(  2.00, 0.00){\line(1,0){4}}
      \put(  2.00, 1.00){\line(1,0){4}}
    \end{picture}$}
  $
  without peak. 
  
  Now let $l<n$ be odd and set $x:=(n-l)/2$ and $p:=2^{.x}.1^{.l}\p n$, then
  the Schensted $Q$-symbol yields a bijection $C_p\to \T^{n,l}$
  (\cite{schensted61}, \cite[Section~4]{schuetzen63}), and the peak
  sets of $\pi$ and $Q(\pi)$ coincide (\cite[Remarque~2]{schuetzen63},
  see also \cite[Theorem~2.1]{foulkes76}).

  Since $\sign(p)=0$, it follows that
  $
  2|\T^{n,l}_{1,e}| = 2|\T^{n,l}_{1,o}|+|C_p\cap \P_0|
  $,
  by Corollary~\ref{comb-cor-1}.
  Denote by $\chi^q$ the irreducible $S_n$-character corresponding
  to $q$ and by $o(q)$ the number of odd columns in the Ferrers diagram of $q$,
  for all $q\p n$, then \cite[Theorem~2.1]{gessel-reutenauer93} and 
  \kref{re-xi-tilde} imply
  $$
  |C_p\cap \P_0|
  =
  \Big(
  \sum_{o(q)=l} \chi^q,\sum_{k=0}^{n-1}\chi^{(n-k).1^{.k}}
  \Big)_{S_n}
  =
  2,
  $$
  since $\sum_{o(q)=l} \chi^q$ is the character of the Lie
  representation of $S_n$ indexed by $p$ (see, for instance,
  \cite[I.8, Example~6(b)]{macdonald95}), and $(n-k).1^{.k}$ is a hook for all
  $k\in\haken{n-1}_0$, hence $\chi^{(n-k).1^{.k}}$ coincides with the
  character of the Foulkes representation of $S_n$ corresponding to
  the composition $1^{.k}.(n-k)$ of $n$ (see \cite{gessel-reutenauer93}).
  This completes the proof.
\end{proof}

 \nopagebreak
\begin{center}
  \textbf{The canonical peak Lie idempotent}
\end{center}

Let $n\in\N$, then
$$
\varrho_n 
:= 
\sum_{q\zerl n}
\frac{(-1)^{\length(q)-1}}{\length(q)}\,\Xi^q
$$
is a Lie idempotent in $\DD_n$ - the \emph{canonical Lie
  idempotent} (\cite{reutenauer85,solomon68}).  By
Proposition~\ref{ex-peak-lie}, it follows that
$$
\tilde{\varrho}_n 
= 
\sum_{q\zerl n}
\frac{(-1)^{\length(q)-1}}{\length(q)}\,\R^q
$$
is (up to the factor $\frac{1}{2}$) a peak Lie idempotent, if $n$
is odd.
Note that
$$
\sum_{n\ge 1}\varrho_n\,t^n
=
\log\Big(1+\sum_{n\ge 1}\Xi^n\,t^n\Big)
\quad
\mbox{ and }
\quad
\sum_{n\ge 1}\tilde{\varrho}_n\,t^n
=
\log\Big(1+\sum_{n\ge 1}\R^n\,t^n\Big),
$$
where $t$ is a variable.

\begin{proposition} \label{canon-1}
  If $n$ is odd, then
  $$
  \tilde{\varrho}_n = 2\sum_{q\zerlodd\,n}
  \frac{(-1)^{(n-\length(q))/2}}{\length(q)}\,\Gamma^q,
  $$
  while $\tilde{\varrho}_n=0$ for even $n$.
\end{proposition}

\begin{proof}
  Standard manipulations based on
  Proposition~\ref{transition-tilde-xi} yield
  \begin{eqnarray*}
    \tilde{\varrho}_n
    & = &
    \sum_{q\zerl n} \frac{(-1)^{\length(q)-1}}{\length(q)}\;\R^q\\
    & = &
    \sum_{q\zerl n} \frac{(-1)^{\length(q)-1}}{\length(q)}\;
    2^{\length(q)}
    \sum_{Q\subseteq\haken{n-1}\backslash(D(q)\cup(D(q)+1))} 
    (-1)^{|Q|}\;\Gamma^Q\\[1mm]
    & = &
    \sum_Q (-1)^{|Q|}\;\Gamma^Q
    \sum_{D\cup (D+1)\subseteq \haken{n}\backslash Q} 
    \frac{(-1)^{|D|}}{|D|+1}\;2^{|D|+1},
  \end{eqnarray*}
  where the sum is taken over all peak sets $Q$ in $\haken{n}$.  But,
  for all such $Q$, setting $E:=Q\cup(Q-1)$ and $m:=n-1-|E|$,
  \begin{eqnarray*}
    \sum_{D\cup (D+1)\subseteq \haken{n}\backslash Q} 
    \frac{(-1)^{|D|}}{|D|+1}\;2^{|D|+1}
    & = &
    \sum_{D\subseteq \haken{n-1}\backslash E} 
    \frac{(-1)^{|D|}}{|D|+1}\;2^{|D|+1}\\
    & = &
    \sum_{k=0}^m
    \binom{m}{k} \frac{(-1)^k}{k+1}\;2^{k+1}\\
    & = &
    \frac{1}{m+1}\,
    \Big(-\sum_{k=1}^{m+1} \binom{m+1}{k} (-2)^k\Big)\\
    & = &
    \frac{1}{m+1}\,\Big(1-(1-2)^{m+1}\Big)\\
    & = &
    \tiptop{\frac{2}{n-|E|}}{\mbox{if $m$ is even,}}{0}{\mbox{otherwise.}}
  \end{eqnarray*}
  But $|E|=2|Q|$ is even, hence $m\equiv n-1$ modulo $2$.  It follows that
  $\tilde{\varrho}_n=0$ for even $n$ and
  $$
  \tilde{\varrho}_n = 2\sum_Q \frac{(-1)^{|Q|}}{n-2|Q|}\,\Gamma^Q,
  $$
  if $n$ is odd, which is another way of stating the claim.
\end{proof}

As a second step, here is an explicit description of the canonical
peak Lie idempotent.

\begin{proposition} \label{peak-canon}
  Let $n\in\N$ be odd, then
  $$
  \textstyle
  \frac{1}{2}\tilde{\varrho}_n 
  = 
  \displaystyle
  \frac{1}{n}\, \sum_{\pi\in S_n}
  (-1)^{\peak(\pi)}\;
  \frac{%
    2\cdot4\cdots (2\cdot\peak(\pi))
  }{%
  (n-2)\cdot(n-4)\cdots(n-(2\cdot\peak(\pi)))}\;\,
  \pi,
  $$
  where $\peak(\pi):=\#\Peak(\pi)$ for all $\pi\in S_n$.
\end{proposition}

\begin{proof}
  By Proposition~\ref{canon-1} and \kref{triangle},
  \begin{eqnarray*}
    \textstyle
    \frac{1}{2}\tilde{\varrho}_n
    & = &
    \sum_{q\zerlodd  n} 
    \frac{(-1)^{(n-\length(q))/2}}{\length(q)}\,\Gamma^q\\
    & = &
    \sum_{q\zerlodd  n} 
    \frac{(-1)^{(n-\length(q))/2}}{\length(q)}
    \sum_{P(q)\subseteq P} \Pi^P\\
    & = &
    \sum_{P} \Big(
    \sum_{\substack{q\zerlodd  n\\[2pt] P(q)\subseteq P}}
    \frac{(-1)^{(n-\length(q))/2}}{\length(q)}
    \Big)\;\Pi^P\\
    & = &
    \sum_{P} \Big( \sum_{Q\subseteq P} \frac{(-1)^{|Q|}}{n-2|Q|} \Big)\;\Pi^P.
  \end{eqnarray*}
  It remains to be shown that, for each peak set
  $P$ in $\haken{n}$ such that $|P|=k$,
  $$
  k_P := \sum_{Q\subseteq P} \frac{(-1)^{|Q|}}{n-2|Q|} =
  \frac{1}{n}\, (-1)^k\frac{2\cdot4\cdots
    (2k)}{(n-2)\cdot(n-4)\cdots(n-2k)}.
  $$
  But $k_P= \sum_{j=0}^k \binom{k}{j} \frac{(-1)^j}{n-2j}=:f_{n,k}$
  depends only on $k$ and $n$, and for these numbers,
  $
  f_{n,0} = \textstyle\frac{1}{n}
  $
  for all $n$, and
  $
  f_{n,k} = f_{n,k-1}-f_{n-2,k-1}
  $
  for all $n\ge 3$, $k\in\N$.  The same recursive formula holds for
  the coefficients in the claim, which thereby follows.
\end{proof}

The canonical peak Lie idempotent will be of importance in
Section~\ref{euler}.

\begin{center}
  \textbf{The Klyachko peak Lie idempotent}
\end{center}

Let $n\in\N$ and $\summe\,D:=\sum_{i\in D} i$ for all
$D\subseteq\haken{n}$, then $ \maj\,\pi := \summe\,\Des(\pi)$
is the \emph{major index} of $\pi$, for all $\pi\in S_n$. Let
$$
\kappa_n(t) := \sum_{\pi\in S_n} t^{\minimaj\pi}\,\pi,
$$
where $t$ is a variable, and denote by $M_z$ the sum of all permutations
$\pi$ in $S_n$ such that $\maj\,\pi\equiv z$ modulo $n$, for all integers $z$.
Once for all, let $\eps$ be a primitive $n$-th root of unity, then due to
Klyachko (\cite{klyachko74}),
$$
\textstyle \frac{1}{n}\kappa_n(\eps) = \displaystyle
\frac{1}{n}\sum_{i=0}^{n-1} \eps^i\,M_i
$$
is a Lie idempotent. As a consequence,
\begin{equation}
  \label{peak-klyachko}
  \textstyle
  \frac{1}{2n}\,\tilde{\kappa}_n(\eps) 
  = 
  \displaystyle
  \frac{1}{2n}\,\sum_{i=0}^{n-1} \eps^i\,(M_i\R^n) 
  =
  \frac{1}{2n}\,\sum_{i=0}^{n-1} \eps^i\,\tilde{M_i}
\end{equation}
is a peak Lie idempotent whenever $n$ is odd. Note that
Proposition~\ref{tilde-delta} implies
$$
\tilde{M}_i
=
\sum_{P}
2^{|P|+1}
\;\#\setofall{D\subseteq\haken{n-1}}{\summe D\equiv_n i,\,
P\subseteq D\sctriangle (D+1)}\;\;
\Pi^P
$$
for all $i\in\haken{n}$, summed over peak sets $P$ in $\haken{n}$.
Unfortunately, we are not able to present a ``cleaner'' description of
the element $\tilde{M}_i$ (or of $\tilde{\kappa}_n(\eps)$) in terms of
the elements $\Pi^Q$.  However, there is a peak analog of
\cite[Proposition~5.1]{LST96}, which seems to be worth mentioning.

If $i\in\haken{n}$ and $d:=\ord\,\eps^i$ denotes the
order of $\eps^i$, then
\begin{equation} \label{LST-wonder}
  \kappa_n(\eps^i)
  =
  \kappa_d(\eps^i)
  \conv
  \cdots
  \conv
  \kappa_d(\eps^i),
\end{equation}
by \cite[Proposition~4.1]{LST96}, hence also
\begin{equation}
  \label{LST-tilde}
  \tilde{\kappa}_n(\eps^i)
  =
  \tilde{\kappa}_d(\eps^i)
  \conv
  \cdots
  \conv
  \tilde{\kappa}_d(\eps^i),
\end{equation}
by Main Theorem~\ref{outer-epi}.

\begin{proposition} \label{simple}
  Let $d,k\in\N$ such that $dk=n$, and let $\beta\in KS_d$ such that
  $\frac{1}{d}\beta$ is a Lie idempotent, then, for
  $\beta^{(k)}
  :=
  \underbrace{\beta\conv\cdots\conv\beta}_{k\;\mbox{\footnotesize factors}}$
  and $\varphi\in\DD_n$,
  $$
  \varphi\beta^{(k)} 
  = 
  c(\varphi)(C_{d^{.k}})\;\beta^{(k)}.
  $$
\end{proposition}

\begin{proof}
  For all $q\zerl n$,
  $$
  \Xi^q\beta^{(k)} 
  = 
  \frac{1}{d^{.k}?}  \Xi^q\omega_{d^{.k}}\beta^{(k)} 
  =
  \frac{1}{d^{.k}?}  \xi^q(C_{d^{.k}})\;\omega_{d^{.k}}\beta^{(k)} 
  =
  \xi^q(C_{d^{.k}})\;\beta^{(k)}\,,
  $$
  by \kref{eq:pr}, \cite[Lemma~1.3]{blelau96} and \cite[(12)]{blelau02}.
  This implies the claim, by linearity.
\end{proof}

Here is the peak analog of \cite[Proposition~5.1]{LST96}.

\begin{corollary} \label{peak-noncyclic}
  Let $n$ be odd, then the linear span $\tilde{\KK}_n$ of the elements
  $\tilde{\kappa}_n(\eps^i)$ ($i\in\haken{n}$) is a left ideal of
  $\DD_n$, and $\tilde{\KK}_n$ is linearly generated by
  $$
  \setofall{\tilde{M}_i}{0\le i\le (n-1)/2}.
  $$
  Furthermore,
  $$
  \tilde{\kappa}_n(\eps^i)\tilde{\kappa}_n(\eps^j) 
  =
  \tiptop{d^{.n/d}?\;\tilde{\kappa_n}(\eps^j)}{\mbox{if
      $d:=\ord\,\eps^j=\ord\,\eps^i$,}} { 0}{\mbox{otherwise,}}
  $$
  and $\tilde{M}_i=\tilde{M}_{n-i}$, for all
  $i,j \in\haken{n}$.  The image of $\tilde{\KK}_n$ under $c$ is
  $$
  \erz{\setofall{\ch_{d^{.k}}}{dk=n}}_K\,.
  $$
\end{corollary}

\begin{proof}
  Let $i,j,d,k,e,l\in\haken{n}$ such that $d=\ord\,\eps^j$,
  $e=\ord\,\eps^i$, and $dk=n=el$, then
  $c(\tilde{\kappa}_n(\eps^i))=\ch_{e^{.l}}$, by \kref{LST-tilde},
  \kref{c-von-lieid} and \kref{chqchr}.  Proposition~\ref{simple} and
  \kref{LST-tilde} thus imply that $\tilde{\KK}_n$ is a left ideal of
  $\DD_n$ as asserted and, in particular,
  $\tilde{\kappa}_n(\eps^i)\tilde{\kappa}_n(\eps^j)
  =\ch_{e^{.l}}(C_{d^{.k}})\,\tilde{\kappa}_n(\eps^j)$.  The multiplication
  rule follows.

  For all $D\subseteq\haken{n-1}$, let $\ol{D}:=\haken{n-1}\backslash D$, then
  on the one hand, Proposition~\ref{multi-tilde-xi} implies
  $$
  \tilde{\Delta}^{D}
  =
  \Delta^{D}\R^n
  =
  \Delta^{D}\Delta^{1^{.n}}\R^n
  =
  \Delta^{\ol{D}}\R^n
  =
  \tilde{\Delta}^{\ol{D}}.
  $$
  On the other hand, for any $\pi\in S_n$,
  $
  \maj\,(\pi\Delta^{1^{.n}})
  = 
  \binom{n}{2}-\maj\,\pi
  \equiv_n
  n-\maj\,\pi
  $,
  since $n$ is odd. It follows that
  $$
  \tilde{M}_i
  = 
  \sum_{\summe D\equiv i} \tilde{\Delta}^D 
  =
  \sum_{\summe D\equiv i} \tilde{\Delta}^{\Dbar} 
  = 
  \sum_{\summe  D\equiv n-i} \tilde{\Delta}^D 
  = 
  \tilde{M}_{n-i}\,.
  $$
  The proof is complete.
\end{proof}

Note that \cite[Proposition~5.1]{LST96} may be obtained along the same
lines, if \kref{LST-wonder} instead of \kref{LST-tilde} is considered.

We conjecture that the elements $\tilde{M}_i$, $0\le i\le (n-1)/2$,
are linearly independent and thus
$$
\dim\,\tilde{\KK}_n = \frac{n+1}{2}.
$$

%
%
%
%
%
%
%
%
%
%
%
%
%
%
%
%
%
%
%
%
%
%
%
%
%

\section{The action on Lie monomials - an algebraic characterization}
\label{alg-approach}

The definition of the descent algebra $\DD_n$ (as well as the definition
of the peak algebra $\PP_n$) is given in purely combinatorial terms.

Due to Garsia and Reutenauer, there is a remarkable result, which
links the descent algebra to the free Lie algebra and the
Poincar\'e-Birkhoff-Witt basis of the free associative algebra
(\cite[Theorem~2.1]{garsia-reutenauer89}).  Based on this result, an
algebraic characterization of $\DD_n$ could be given (ibid.,
Theorem~4.5). The aim of this section is to derive an analogous
algebraic characterization of the peak algebra.

Recall that
$
L(X) = \bigoplus_{n\ge 0} L_n(X)
$
denotes the free Lie algebra over the alphabet $X$, contained in the free
associative algebra $A(X)$. If $n\in\N$, $q=q_1\wldots q_k\zerl n$ and
$P_1\in L_{q_1}(X),\ldots,P_k\in L_{q_k}(X)$, then, for all
$S=\{i_1,\ldots,i_m\}\subseteq \haken{k}$ such that $i_1<\cdots <i_m$, put
$$
q_S := q_{i_1}\wldots q_{i_k} \quad\mbox{ and }\quad P_S := P_{i_1}\cdots P_{i_k}\,.
$$
Furthermore, denote by $\SP(\haken{k})$ the set of all set
partitions $(S_1,\ldots,S_l)$ of $\haken{k}$.  With these notations,
\cite[Theorem~2.1]{garsia-reutenauer89} reads as follows:

\begin{theorem} \label{gr2.1}
  Let $n\in\N$, $q=q_1\wldots q_k,r=r_1\wldots r_l\zerl n$ and $P_i\in
  L_{q_i}(X)$ for all $i\in\haken{k}$, then
  $$
  \Xi^rP_1\cdots P_k 
  = 
  \sum_{%
    \substack{(S_1,\ldots, S_l)\in\SP(\haken{k})\\[2pt]  q_{S_i}\zerl r_i}}
  P_{S_1}\cdots P_{S_l}\,,
  $$
  via Polya action.
\end{theorem}

As a consequence, for all $\varphi\in KS_n$, it follows that
$\varphi\in\DD_n$ if and only if
$$
\varphi P_1\cdots P_k \in \Erz{\setofall{P_{1\pi}\cdots
    P_{k\pi}}{\pi\in S_k}}_K\,,
$$
for all $q=q_1\wldots q_k\zerl n$, $P_1\in L_{q_1}(X),\ldots,P_k\in
L_{q_k}(X)$ (see \cite[Theorem~4.5]{garsia-reutenauer89}).  In
particular, $\DD_n$ acts on the linear span of the elements
$P_{1\pi}\cdots P_{k\pi}$, $\pi\in S_k$, from the left.  One is
tempted to ask if $\PP_n$ may be characterized by \emph{the way} the
elements of $\PP_n$ act on these linear spaces; the answer is yes. The
surprising characterization is again related to parity and based on
the following application of Theorem~\ref{gr2.1}.

\begin{lemma} \label{gerade-null}
  Let $n\in\N$, $q=q_1\wldots q_k\zerl n$ and $P_i\in L_{q_i}(X)$ for
  all $i\in\haken{k}$.  If $q_1$ is even, then
  $$
  \tilde{\Xi}^nP_1\cdots P_k=0.
  $$
\end{lemma}

\begin{proof}
  Denoting by $\SP(\haken{k})_{odd}$ the set of all
  $(S_1,\ldots,S_m)\in\SP(\haken{k})$ such that $\summe q_{S_m}$ is
  odd, it follows from Proposition~\ref{transition-tilde-xi} and
  Theorem~\ref{gr2.1} that
  \begin{eqnarray*}
    \tilde{\Xi}^nP_1\cdots P_k
    & = &
    2\sum_{%
      \substack{s\zerl n\\[2pt] s^\dagger\;\mbox{\tiny odd}}} 
    (-1)^{n-\length(s)}\;\Xi^sP_1\cdots P_k\\
    & = &
    2\sum_{%
      \substack{s\zerl n\\[2pt] s^\dagger\;\mbox{\tiny odd}}}
    (-1)^{n-\length(s)}
    \sum_{%
      \substack{(S_1,\ldots, S_{\length(s)})\in\SP(\haken{k})\\[2pt]
        q_{S_i}\zerl s_i}}
    P_{S_1}\cdots P_{S_{\length(s)}}\\
    & = &
    2  \sum_{(S_1,\ldots, S_m)\in\SP(\haken{k})_{odd}} 
    (-1)^{n-m} P_{S_1}\cdots P_{S_m}\,.\\
  \end{eqnarray*}
  Let $A$ be the set of all $(S_1,\ldots, S_m)\in\SP(\haken{k})_{odd}$
  such that $S_i=\{1\}$ for an index $i\in\haken{m}$, and set
  $B:=\SP(\haken{k})_{odd}\backslash A$.  Note that $S_i=\{1\}$
  implies that $i<m$, for all $(S_1,\ldots,
  S_m)\in\SP(\haken{k})_{odd}$, since $q_1$ is even and $\summe
  q_{S_m}$ is odd.  For all $S=(S_1,\ldots, S_m)\in A$, choose
  $i\in\haken{m-1}$ such that $S_i=\{1\}$ and define
  $$
  \tilde{S} 
  := 
  (S_1,\ldots,S_{i-1},S_i\cup S_{i+1},S_{i+2},\ldots,S_m) 
  \in 
  B.
  $$
  Conversely, for all $T=(T_1,\ldots, T_m)\in B$, let $i\in\haken{m}$
  such that $1\in T_i$ and define
  $$
  \hat{T} 
  :=
  (T_1,\ldots,T_{i-1},\{1\},T_i\backslash\{1\},T_{i+1},\ldots,T_m) 
  \in
  A.
  $$
  The mapping $A\to B,\,S\mapsto \tilde{S}$ is a bijection with
  inverse given by $T\mapsto \hat{T}$.  Furthermore, for all
  $S=(S_1,\ldots, S_m)\in A$, if $\tilde{S}=(\tilde{S}_1,\ldots,
  \tilde{S}_{m-1})$, then
  $$
  P_{S_1}\cdots P_{S_m} = P_{\tilde{S}_1}\cdots P_{\tilde{S}_{m-1}}\,,
  $$
  hence
  \begin{eqnarray*}
    \tilde{\Xi}^nP_1\cdots P_k
    & = &
    2  \sum_{(S_1,\ldots, S_m)\in\SP(\haken{k})_{odd}} 
    (-1)^{n-m} P_{S_1}\cdots P_{S_m}\\
    & = &
    2  \sum_{(S_1,\ldots, S_m)\in A} 
    (-1)^{n-m}
    \Big(
     P_{S_1}\cdots P_{S_m}- P_{\tilde{S}_1}\cdots P_{\tilde{S}_{m-1}}
    \Big)\\
    & = &
    0
  \end{eqnarray*}
  as asserted.
\end{proof}

We proceed by transferring Theorem~\ref{gr2.1} to the group ring $KS_n$.  The
symmetric group $S_n$ acts on $KS_n$ by left multiplication, and on
$A_n(X)$ via Polya action.  Mapping $\pi\in S_n$ onto the multi-linear
word
$
(1\pi).(2\pi)\wldots(n\pi)\in\N^*
$
extends to an $S_n$-left module monomorphism $\iota_n:KS_n\to A_n(\N)$, by
linearity. For any $\pi\in S_n$, $\sigma\in S_m$, define 
$\pi\#\sigma\in S_{n+m}$ by
$$
i(\pi\#\sigma)
:=
\tiptop{i\pi}{\mbox{if $i\le n$,}}{(i-n)\sigma+n}{\mbox{if $i> n$,}}
$$
for all $i\in\haken{n+m}$.
Bilinear extension gives an associative product $\#$ on $KS$.
Let $\gamma$ be a Lie series in $\DD$ and put
$
\gamma^q
:=
\gamma_{q_1}\#\cdots\#\gamma_{q_k}
$
for all $q=q_1\wldots q_k\in \N^*$, then
\begin{equation}
  \label{up-down}
  \gamma_q=\gamma^q\Xi^q
\end{equation}
(see, for instance, \cite[Proposition~2.1]{schocker-lia}).
Furthermore, there is a Lie monomial $P_i\in L_{q_i}(X)$ for each
$i\in\haken{k}$ as in Theorem~\ref{gr2.1} such that
\begin{equation}
  \label{transfer}
  \gamma^q\iota_n
  =
  P_1\cdots P_k\,.
\end{equation}
As a consequence, Theorem~\ref{gr2.1} implies that 
$\gamma_r\gamma_q=\gamma^r(\Xi^r\gamma^q)\Xi^q=0$
unless there is composition $\dot{q}\zerl n$ such that $q\ass \dot{q}\zerl r$.
In this case, write $q\assozerl r$, that is,
\begin{equation}
  \label{eq:lia}
  \gamma_r\gamma_q\neq 0\mbox{ implies }q\assozerl r.
\end{equation}

With these preparations, the preceding lemma yields an \emph{internal}
characterization of the elements $\varphi\in\PP_n$ --- based on the
action of $\varphi$ on $\DD_n$ itself from the left.

\begin{corollary} \label{internal-char}
  Let $n\in\N$, $\varphi\in\DD_n$ and $\gamma$ be a Lie series in
  $\DD$, then $\varphi\in\PP_n$ if and only if
  $$
  \varphi\gamma_q=0
  $$
  for all $q=q_1\wldots q_l\zerl n$ such that $q_1$ is even.
\end{corollary}

\begin{proof}
  For any $\varphi\in\PP_n$, there exists an $\alpha\in\DD_n$ such
  that $\varphi=\alpha\tilde{\Xi}^n$, by Main Theorem~\ref{outer-epi},
  hence \kref{up-down}, \kref{transfer} and Lemma~\ref{gerade-null}
  imply $ \varphi\gamma_q = \alpha(\tilde{\Xi}^n\gamma^q)\Xi^q = 0$
  whenever $q_1$ is even.
  
  Conversely, let $\varphi\in\DD_n$ such that $\varphi\gamma_q=0$
  whenever $q_1$ is even.  If $\beta$ is another Lie series in $\DD$,
  it follows that $ \varphi\beta_q =
  \textstyle\frac{1}{q?}\varphi\gamma_q\beta_q = 0 $ for all $q\zerl
  n$ such that $q_1$ is even, by \kref{eq:pr}.  Without loss of
  generality, we may thus assume that $\gamma_n\in\PP_n$, for all odd
  $n$, hence $\PP_n=\erz{\setofall{\gamma_q}{q\zerlodd n}}_K$, by Main
  Theorem~\ref{high-lie-id-basis}.
  
  Let $\varphi=\sum_{r\zerl n} a_r\gamma_r$ and assume that
  $\varphi\notin\PP_n$, then there is a composition $s$ of $n$ such
  that $a_s\neq 0$ and $s$ is not odd. 
  Choose $s$ of minimal length
  and $t\ass s$ such that $t_1$ is even, then
  $$
  0
  =
  \varphi \gamma_t
  =
  \sum_{r\zerlodd  n}a_r\gamma_r\gamma_t
  +
  \sum_{%
    \substack{r\zerl n\\[2pt] r\;\,\mbox{\tiny not odd}}}
  a_r\gamma_r\gamma_t\\
  =
  \sum_{%
    \substack{r\zerl n\\[2pt] r\;\,\mbox{\tiny not odd}}}
  a_r\gamma_r\gamma_t\,,
  $$
  by the part already proven.  But $a_r\gamma_r\gamma_t\neq 0$
  implies $a_r\neq 0$ and $s\ass t\assozerl r$, for all $r\zerl n$, by
  \kref{eq:lia}.  If, additionally, $r$ is not odd, then
  $\length(s)\le\length(r)$, due to the choice of $s$, hence
  $\length(s)=\length(r)$ and $s\ass r$.
  It follows that 
  $$
  \sum_{%
    \substack{r\zerl n\\[2pt] r\;\,\mbox{\tiny not odd}}}
  a_r\gamma_r\gamma_t
  =
  r? \sum_{s\ass r} a_r\gamma_r\,,
  $$
  But the elements $\gamma_r$, $r\ass s$, are linearly independent,
  which implies $a_r=0$ for all $r\ass s$, a contradiction.  This
  shows $\varphi\in\PP_n$.
\end{proof}

Combining Lemma~\ref{gerade-null} and Corollary~\ref{internal-char}
leads to the following algebraic characterization of the left ideal
$\PP_n$ of $\DD_n$.

\begin{mtheorem}
  Let $n\in\N$ and $\varphi\in\DD_n$, then 
  $\varphi\in\PP_n$ if and only if
  $$
  \varphi\,P_1\cdots P_k=0
  $$
  for all $q=q_1\wldots q_k\zerl n$ such that $q_1$ is even and all
  $P_1\in L_{q_1}(X),\ldots,P_k\in L_{q_k}(X)$.
\end{mtheorem}

\begin{proof}
  The necessity part follows from Lemma~\ref{gerade-null} and Main
  Theorem~\ref{outer-epi}.  Conversely, $ \varphi\,P_1\cdots P_k=0 $
  for all the Lie monomials $P_1,\ldots,P_k$ mentioned in the claim
  implies $\varphi\gamma^q=0$ for all $q\zerl n$ such that $q_1$ is
  even, by \kref{transfer}, hence also
  $\varphi\gamma_q=\varphi\gamma^q\Xi^q=0$.  Now $\varphi\in\PP_n$
  follows from Corollary~\ref{internal-char}.
\end{proof}

It seems to be worth mentioning that putting together the preceding theorem and
the Garsia-Reutenauer characterization of $\DD_n$ yields:

Any $\varphi\in KS_n$ is constant on each peak class in $S_n$ if and only if,
via Polya action, $\varphi$ acts on the linear span
$
\erz{\setofall{P_{1\pi}\cdots P_{k\pi}}{\pi\in S_k}}_K
$
for all $q=q_1\wldots q_k\zerl n$, $P_1\in L_{q_1}(X),\ldots,P_k\in
L_{q_k}(X)$, and annihilates $P_1\cdots P_k$ whenever $q_1$ is even.

%
%
%
%
%
%
%
%
%
%
%
%
%
%
%
%
%
%
%
%
%
%
%
%
%

\section{The Eulerian peak algebra} \label{euler}

Let $n\in\N$ and consider the sums of all permutations in $S_n$ with a
given number of descents:
$$
\Delta^{n,k} 
:= 
\sum_{\substack{q\zerl n\\[2pt]\length(q)=k}} \Delta^q
\qquad(k\in\haken{n}).
$$
The linear span $\EE_n$ of these elements is an $n$-dimensional
commutative sub-algebra of $\DD_n$ (see, for instance,
\cite{cellini95,garsia-reutenauer89}), usually referred to as the
\emph{Eulerian sub-algebra} of $\DD_n$.  The aim of this section is to
study the peak variant of $\EE_n$. Recall first that a second basis of
$\EE_n$ is constituted by the elements
$$
\Xi^{n,k} 
:= 
\sum_{\substack{q\zerl n\\[2pt]\length(q)=k}} \Xi^q
\qquad(k\in\haken{n}),
$$
and that the canonical Lie idempotent $\varrho_n$ has been defined
in Section~\ref{PLI}.  If we set
$\varrho_q:=\varrho_{q_1}\conv\cdots\conv\varrho_{q_k}$, for all
$q=q_1\wldots q_k\in\N^*$, then due to Reutenauer
(\cite{reutenauer85}, see also \cite{garsia89}), the elements
$$
\varrho_{n,k} 
:= 
\frac{1}{k!}\sum_{\substack{q\zerl n\\[2pt]\length(q)=k}} \varrho_q
\qquad(k\in\haken{n})
$$
are mutually orthogonal idempotents in $\EE_n$, and
$\sum_{k=1}^n\varrho_{n,k}=\id_n$.

Now consider the sums of all permutations in $S_n$ with a
given number of peaks:
$$
\Pi^{n,k} 
:= 
\sum_{%
\substack{P\;\mbox{\footnotesize peak set in }\haken{n}\\[2pt] |P|=k}} 
\Pi^P,
$$
and denote by $\tilde{\EE}_n$ the linear span of these elements.
Furthermore, 
let
$$
\Gamma^{n,k} 
:= 
\sum_{%
  \substack{P\;\mbox{\footnotesize peak set in }\haken{n}\\[2pt] |P|=k}} 
\Gamma^P
$$
for all $k\in\N_0$.

\begin{proposition} \label{euler-transition}
  The sets
  $$
  \setofall{\Pi^{n,k}}{0\le k\le (n-1)/2} 
  \quad\mbox{ and }\quad
  \setofall{\Gamma^{n,k}}{0\le k\le (n-1)/2}
  $$
  are linear bases of $\tilde{\EE}_n$, and
  \begin{equation}
    \label{eq:euler-transition}
    \Gamma^{n,k} = \sum_{j\ge k} \textstyle\binom{j}{k}\, \Pi^{n,j}
    \quad\mbox{ and }\quad 
    \displaystyle
    \Pi^{n,k} = \sum_{j\ge k} (-1)^{j-k}\textstyle\binom{j}{k}\, \Gamma^{n,j}
  \end{equation}
  for all $0\le k\le (n-1)/2$.
  Furthermore, 
  $
  \R^{n,k}
  =
  \Xi^{n,k}\R^n
  \in
  \tilde{\EE}_n
  $
  for all $k\in\haken{n}$.  
\end{proposition}

\begin{proof}
  By definition, $\setofall{\Pi^{n,k}}{0\le k\le (n-1)/2}$ is a linear
  basis of $\tilde{\EE}_n$.  
  
  The identities \kref{eq:euler-transition} are immediate from
  \kref{triangle} and \kref{triangle-back}.  In particular, the
  elements $\Gamma^{n,k}$, $0\le k\le (n-1)/2$ are also linearly
  independent.  
  Finally,
  $$
  \R^{n,k}
  =
  \sum_{\substack{q\zerl n\\[2pt]\length(q)=k}} \R^q\\
  = 
  2^k \sum_{j\ge 0}  (-1)^j \binom{n-1-2j}{k-1} \;\Gamma^{n,j}
  \in
  \tilde{\EE}_n
  $$
  is readily seen from \kref{xi-tilde} and
  Proposition~\ref{transition-tilde-xi}.
\end{proof}

Now observe that
\begin{equation}
  \label{def-rhonk}
  \tilde{\varrho}_{n,k}
  =
  \varrho_{n,k}\R^n
  =
  \frac{1}{k!}
  \sum_{\substack{q\zerl\, n\\[2pt]\length(q)=k}} \tilde{\varrho}_q
  =
  \frac{1}{k!}
  \sum_{\substack{q\zerlodd\, n\\[2pt]\length(q)=k}} \tilde{\varrho}_q\,,
\end{equation}
for all $k\in\haken{n}$, where the latter equality is a consequence of
Proposition~\ref{canon-1}.  In particular, $\tilde{\varrho}_{n,k}\neq
0$ if and only if $k\equiv n$ modulo $2$, since $q\zerlodd n$ implies
$\length(q)\equiv n$ modulo $2$.

\begin{mtheorem} \label{euler-alg}
  $\tilde{\EE}_n$
  is a commutative sub-algebra of $\PP_n$ of dimension $(n+1)/2$, if
  $n$ is odd, and $n/2$, if $n$ is even.  
  As an algebra, $\tilde{\EE}_n$ is generated by the single element 
  $$
  \Pi^{n,0}
  =
  \Pi^\emptyset
  =
  \sum_{k=0}^n \Delta^{1^{.k}.(n-k)}.
  $$
  Furthermore, 
  $
  \setofall{\frac{1}{2^k}\tilde{\varrho}_{n,k}}{%
    k\in\haken{n},\,k\equiv n\mbox{ modulo }2}
  $
  is a linear basis of $\tilde{\EE}_n$ consisting of mutually
  orthogonal idempotents. In particular, 
  $c|_{\tilde{\EE}_n}$ is a monomorphism of algebras,
  mapping $\tilde{\EE}_n$ onto the linear span of the elements
  $$
  \cha_{n,k}
  :=
  \sum_{\substack{p\podd\, n\\[2pt]\length(p)=k}} \cha_p
  \qquad(k\equiv n\mbox{ modulo }2).
  $$
\end{mtheorem}

The first statement in the preceding theorem is due to Doyle and Rockmore
(\cite{doyle-rockmore02}).

\begin{proof}
  By \kref{c-von-lieid}, Main Theorem~\ref{p-c-im} and \kref{chqchr},
  $c(\varrho_q)=\frac{1}{q_1\cdots q_k}\,\ch_q$
  for all $q=q_1\wldots q_k\zerl n$.
  For each $k\in\haken{n}$ and $p\p n$, it follows that
  $
  c(\varrho_{n,k})(C_p)
  =
  \frac{1}{k!}\,\sum_{\length(q)=k}c(\varrho_q)(C_p)
  $
  does not vanish if and only if $\length(p)=k$, which implies
  $$
  c(\varrho_{n,k})
  =
  \sum_{\substack{p\p n\\[2pt]\length(p)=k}} \cha_p\,,
  $$
  since the elements $c(\varrho_{n,k})$ are pairwise orthogonal idempotents
  in $\Cl_K(S_n)$ summing up to 
  $c(\id_n)=c(\Xi^n)=\xi^n=\sum_{p\p n}\cha_p$,
  by Theorem~\ref{sol-rad}.
  
  As another consequence of Theorem~\ref{sol-rad} and Main
  Theorem~\ref{p-c-im},
  $$
  \textstyle
  c(\frac{1}{2^k}\tilde{\varrho}_{n,k})
  =
  \frac{1}{2^k}c(\varrho_{n,k})c(\R^n)
  =
  \frac{1}{2^k}
  \displaystyle
  \sum_{\substack{p\p n\\[2pt]\length(p)=k}} \cha_p
  \sum_{p\podd  n} 2^{\length(p)}\,\cha_p
  =
  \cha_{n,k}
  $$
  as asserted. Furthermore, $\varrho_{n,k}\in\EE_n$, hence
  $\tilde{\varrho}_{n,k}$ is a linear combination of the elements
  $\R^{n,k}$.  Proposition~\ref{euler-transition} implies
  $\tilde{\varrho}_{n,k}\in\tilde{\EE}_n$ and, in particular,
  $$
  \dim\,c(\tilde{\EE}_n)
  \ge
  \dim\,\erz{\setofall{\cha_{n,k}}{k\equiv n\mbox{ modulo }2}}_K
  =
  \dim\,\tilde{\EE}_n\,.
  $$
  In other words, $c|_{\tilde{\EE}_n}$ is injective, mapping
  $\tilde{\EE}_n$ onto the linear span of the elements $\cha_{n,k}$,
  $k\equiv n$ modulo $2$.  But these elements are mutually orthogonal
  idempotents, hence so are the elements
  $\frac{1}{2^k}\tilde{\varrho}_{n,k}$, $k\equiv n$ modulo $2$, by
  Main Theorem~\ref{p-c-im}. In particular, $\tilde{\EE}_n$ is a
  commutative sub-algebra of $\PP_n$.

  Let $\ZZ_n$ be the sub-algebra of $\tilde{\EE}_n$ generated by
  $\Pi^\emptyset=\frac{1}{2}\R^n$,
  then $c(\ZZ_n)$ is a sub-algebra of $\Cl_K(S_n)$ containing
  $$
  c(\R^n)
  =
  \sum_{p\podd  n} 2^{\length(p)}\,\cha_p
  =
  \sum_{k\equiv n} 2^k\,\cha_{n,k}\,.
  $$
  The transition matrix expressing the powers of $c(\R^n)$ in
  $\Cl_K(S_n)$ in the elements $\cha_{n,k}$ is thus a Vandermonde
  matrix, which particularly implies
  $
  \dim\,\ZZ_n
  =
  \dim\, c(\ZZ_n)
  =
  \dim\,c(\tilde{\EE}_n)
  =
  \dim\,\tilde{\EE}_n
  $,
  hence $\ZZ_n=\tilde{\EE}_n$. This completes the proof.
\end{proof}

\begin{corollary}
  The element 
  $e:=\sum_{k=1}^n \frac{1}{2^k}\tilde{\varrho}_{n,k}$
  is an idempotent in $\PP_n$ such that $\PP_n=\DD_ne$.
\end{corollary}

\begin{proof}
  The element $e$ is an idempotent, by \kref{def-rhonk} and Main
  Theorem~\ref{euler-alg}, and
  $
  \DD_ne
  =
  \DD_n\,(\sum_{k=1}^n\varrho_{n,k})\,\R^n
  =
  \DD_n\R^n
  =
  \PP_n\,,
  $
  by Main Theorem~\ref{outer-epi}.
\end{proof}

\begin{remark}
  As in the case of $\PP_n$ itself (Section~\ref{comb-approach}),
  there is a direct combinatorial approach to the essential part of
  Main Theorem~\ref{euler-alg}, as will be briefly explained now.
  
  Let $\pi,\sigma\in S_n$, then $\sigma$ is called a \emph{peak number
    neighbor} of $\pi$ if $\sigma$ is a peak neighbor of $\pi$ or
  $\sigma=\pi\tau_{n,n-1}$ and $\{(n-1)\pi^{-1},n\pi^{-1}\}\subseteq
  \{2,\ldots,n-1\}$. In this case, write $\sigma\peaknumn\pi$.
  Furthermore, denote by $\peaknumr$ the smallest equivalence relation
  on $S_n$ containing $\peaknumn$.
  
  In analogy to Lemma~\ref{peak-rel}, there is the following result:
  \begin{center}
    \emph{
      $\peak(\pi)=\peak(\sigma)$ if and only
      if $\pi\peaknumr \sigma$, for all $\pi,\sigma\in S_n\,$.}
  \end{center}
  Set
  $
  F(i,j,\pi)
  :=
  \setofall{(\alpha,\beta)\in S_n\times S_n}
  {\peak(\alpha)=i,\,\peak(\beta)=j,\,\alpha\beta=\pi}
  $,
  for all $i,j\in\N_0$ and observe that
  $\tilde{\EE}_n$ is multiplicatively closed if
  $\# F(i,j,\pi)=\# F(i,j,\sigma)$
  for all $\pi,\sigma\in S_n$ such that $\peak(\pi)=\peak(\sigma)$.
  The latter condition is equivalent to $\pi\peaknumr\sigma$,
  as was mentioned above, and it suffices to consider the case
  where $\pi\peaknumn\sigma$.
  
  If $\pi\peakn\sigma$, then $\Peak(\pi)=\Peak(\sigma)$, and $\#
  F(i,j,\pi)=\# F(i,j,\sigma)$ follows from the fact that $\PP_n$ is a
  sub-algebra of $\DD_n$.
  
  The case where $\sigma=\pi\tau_{n,n-1}$ and
  $\{(n-1)\pi^{-1},n\pi^{-1}\}\subseteq \{2,\ldots,n-1\}$ remains.  But,
  in this case, $ \iota_\pi:F(i,j,\pi)\lra F(i,j,\sigma) $, defined by
  $$
  (\alpha,\beta)
  \lmt
  \tiptop{(\alpha,\beta\tau_{n,n-1})}{%
    \mbox{if $\beta\tau_{n,n-1}\peaknumn\beta$,}}
  {(\alpha(\beta\tau_{n,n-1}\beta^{-1}),\beta)}{\mbox{otherwise,}}
  $$
  is a bijection (with inverse given by $\iota_\sigma$).
  
  Thus it follows again that $\tilde{\EE}_n$ is a sub-algebra of $\PP_n$.
\end{remark}

%
%
%
%
%
%
%
%
%
%
%
%
%
%
%
%
%
%
%
%
%
%
%
%
%

\section{Principal indecomposables and the Jacobson radical} \label{radical}

The peak algebra $\PP_n$ is not only a sub-module of the regular
$\DD_n$-left module, but even a direct summand (Main
Theorem~\ref{dir-summ}).  Basic information about the structure of
$\PP_n$ is thus easily deduced from the corresponding results on
$\DD_n$.  This includes a decomposition into indecomposable
$\PP_n$-left modules and a description of the Jacobson radical of
$\PP_n$ (Corollary~\ref{rad-lambda}) as well as of its Cartan
invariants (Corollary~\ref{cartan-pn}).  More effort is required to
collect informations on the descending Loewy series of $\PP_n$, which
allows to determine the nil-potency index of the Jacobson radical of
$\PP_n$ (Theorem~\ref{nilindex}).

The results on the structure of $\DD_n$ derived in \cite{blelau96} and
\cite{blelau02} are of crucial importance for what follows.
Furthermore, we should like to point out again that $\PP_n$ is an
algebra \emph{without} identity; standard arguments from
representation theory may thus be applied with care only.

Throughout this section, $\gamma$ is a Lie series in $\DD$.  Recall
that $\gamma$ extends to an isomorphism of algebras
$
\gamma:A(\N)\to (\DD,\conv),\,\sum a_q q\mapsto \sum a_q\gamma_q
$
(see \kref{gamma-iso}).
To start with, here is an immediate consequence of 
\cite[Lemma~1.3]{blelau96}.

\begin{proposition} \label{indec-dn}
  Let $n\in\N_0$, then
  $\gamma$ maps each rearrangement class
  $\setofall{q}{q\ass p}$ ($p\p n$) onto
  the indecomposable $\DD_n$-left module
  $$
  \Lambda^p 
  := 
  \DD_n\gamma_p
  =
  \langle\,\setofall{\gamma_{q}}{q\ass p}\,\rangle_K,
  $$
  and
  $
  \DD_n = \bigoplus_{p\p n} \Lambda^p
  $.
  Furthermore, $\Lambda^p\cong \DD_n\omega_p$ as a $\DD_n$-left
  module.
\end{proposition}

\begin{proof}
  Let $p\p n$, then the mapping
  $
  \DD_n\gamma_p\lra\DD_n\omega_p,\, 
  \varphi
  \lmt
  \textstyle\frac{1}{p?}\varphi\omega_p
  $
  is an isomorphism of $\DD_n$-left modules with inverse
  $\varphi\lmt \frac{1}{p?}\varphi\gamma_p$, by \kref{eq:pr}.  In
  particular, $\Lambda^p$ is an indecomposable $\DD_n$-left module,
  since
  $$
  \DD_n\omega_p 
  = 
  \langle\,\setofall{\omega_q}{q\ass p}\,\rangle_K
  $$
  is indecomposable (\cite[Lemma~1.3]{blelau96}),
  and
  $
  \Lambda^p 
  = 
  \DD_n\omega_p\gamma_p 
  =
  \langle\,\setofall{\gamma_q}{q\ass p}\,\rangle_K
  $,
  by \kref{eq:pr} again.  Finally, 
  $
  \DD_n 
  = 
  \langle\,\setofall{\gamma_q}{q\zerl n}\,\rangle_K
  =
  \bigoplus_{p\p n} \Lambda^p
  $,
  as $\gamma:A(\N)\to\DD$ is an isomorphism.
\end{proof}

As a consequence of \kref{c-von-lieid}, Main Theorem~\ref{p-c-im}
and \kref{chqchr},
\begin{equation}
  \label{c-von-highlieid}
  c(\gamma_q)
  =
  \ch_q
\end{equation}
for all $q\in\N^*$, which leads to

\begin{proposition} \label{lambda-q-regu}
  Let $n\in\N_0$ and $p\p n$, then $\Lambda^p$ contains a unique maximal
  $\Lambda^p$-left sub-module, namely
  $
  \Lambda^p\cap\ker\,c 
  =
  \erz{\setofall{\gamma_p-\gamma_q}{q\ass p}}_K
  $.

  In particular, $\Lambda^p$ is indecomposable as a $\Lambda^p$-left module.
\end{proposition}

\begin{proof}
  The stated equality is immediate from Proposition~\ref{indec-dn} and
  \kref{c-von-highlieid}, hence $K:=\Lambda^p\cap\ker\,c$ is a
  $\Lambda^p$-sub-module of $\Lambda^p$ of co-dimension $1$ and thus
  maximal.
  
  Let $M$ be an arbitrary maximal $\Lambda^p$-sub-module of
  $\Lambda^p$.  If $m=\sum_{q\ass p} k_q\gamma_q\in M$ and 
  $s\ass p$, then
  $$
  p?\,\Big(\sum_{q\ass p} k_q\Big)\,\gamma_s = \gamma_sm \in M,
  $$
  by \kref{eq:pr}, hence $\sum_{q\ass p} k_q=0$, since $M\neq
  \Lambda^p$.  It follows that $c(m)=(\sum_{q\ass p} k_q)\ch_p=0$, by
  \kref{c-von-highlieid}, thus $M\subseteq K$ and even $M=K$.
\end{proof}

From now on, assume that $\gamma_n\in\PP_n$, for all odd $n$.
Furthermore, $n\in\N$ is fixed, and the set of all odd partitions of
$n$ is denoted by $\Part_{odd}(n)$.

\begin{mtheorem} \label{dir-summ}
  $\PP_n=\bigoplus_{p\podd\, n} \Lambda^p$ is a decomposition of
  $\PP_n$ into indecomposable $\PP_n$-left modules.  In particular,
  $\PP_n$ is a direct summand of the regular $\DD_n$-left module:
  $$
  \DD_n = \bigoplus_{p\p n} \Lambda^p = \PP_n\oplus\bigoplus_p
  \Lambda^p,
  $$
  where the latter sum runs over all
  $p\in\Part(n)\backslash\Part_{odd}(n)$.
\end{mtheorem}

\begin{proof}
  By Main Theorem~\ref{high-lie-id-basis} and
  Proposition~\ref{indec-dn}, indecomposability of $\Lambda^p$ as a
  $\PP_n$-left module remains, for $p\podd n$.  But this follows
  directly from Proposition~\ref{lambda-q-regu}.
\end{proof}

For any finite-dimensional $K$-algebra $A$ and any $A$-left module
$M$, the \emph{$A$-radical}
$
\Rad_AM
$
of $M$ is the intersection of all maximal $A$-sub-modules of $M$.
In particular, 
$
\Rad\,A=\Rad_AA
$
is the \emph{Jacobson radical} of $A$.

\begin{corollary} \label{rad-lambda}
  $
  \Rad_{\PP_n}\Lambda^p=\ker\,c\cap\Lambda ^p=\Rad_{\DD_n}\Lambda^p
  $,
  for all $p\podd  n$.  In particular,
  $$
  \Rad\,\PP_n 
  = 
  \bigoplus_{p\podd n}\Rad_{\PP_n}\Lambda^p 
  =
  \ker\,c|_{\PP_n}
  =
  \erz{\setofall{\gamma_q-\gamma_p}{q\ass p\podd n}}_K\,,
  $$
  and $\PP_n/\Rad\,\PP_n$ is commutative and of dimension
  $\#\Part_{odd}(n)$.
\end{corollary}

\begin{proof}
  The first claim follows from Propositions~\ref{indec-dn},
  \ref{lambda-q-regu} and \cite[Corollary~1.6]{blelau96}.  
  As a consequence,
  $$
  \Rad\,\PP_n 
  \subseteq 
  \bigoplus_{p\podd n}\Rad_{\PP_n}\Lambda^p 
  =
  \ker\,c|_{\PP_n}
  =
  \erz{\setofall{\gamma_q-\gamma_p}{q\ass p\podd n}}_K\,,
  $$
  and $\gamma_p\Rad_{\PP_n}\Lambda^p=0$, by \kref{eq:pr}.

  Let $M$ be a maximal sub-module of $\PP_n$ and $p\podd n$.
  Assume that $\Rad_{\PP_n}\Lambda^p$ is not contained in
  $M$, then $M+\Rad_{\PP_n}\Lambda^p=\PP_n$.
  In particular, there are $m\in M$ and $r\in\Rad_{\PP_n}\Lambda^p$
  such that $\gamma_p=m+r$, which implies
  $$
  \textstyle
  p?\,\gamma_p
  =
  \gamma_p(m+r)
  =
  \gamma_pm
  \in
  M,
  $$
  hence even $\Lambda^p\subseteq M$, by \kref{eq:pr} --- a
  contradiction.  This shows the remaining inclusion.
  
  An application of Main Theorem~\ref{p-c-im} completes the proof.
\end{proof}

Recall from \cite[p.~319]{blelau02} that
$$ 
M^p := \Lambda^p/\Rad_{\DD_n}\Lambda^p 
$$
is a one-dimensional irreducible $\DD_n$-module, for all $p\p n$,
and that the modules $M^p$, $p\p n$, represent the isomorphism classes
of irreducible $\DD_n$-modules.

If $n=1$, then $\PP_n=\DD_n$. For the remainder, assume that
$n>1$, then, for instance,
$$
M^0:=M^{2.1^{.(n-2)}}
$$
is a trivial $\PP_n$-module, reflecting the fact that $\PP_n$ is an
algebra without identity.  More generally, if $V$ is a $\PP_n$-module,
then
$$
\ol{V} := \setofall{v\in V}{\varphi\,v=0\mbox{ for all
    }\varphi\in\PP_n}
$$
is a zero sub-module of $V$, that is
$
\ol{V} \cong_{\PP_n} (\dim\,\ol{V}) M^0
$.

\begin{corollary}
  $M^p$, $p\podd  n$, and $M^0$ represent the isomorphism classes of
  irreducible $\PP_n$-modules.  Furthermore,
  $$
  M^0 \cong M^q
  $$
  (as $\PP_n$-left modules) for all $q\p n$ containing an even
  letter.
\end{corollary}

\begin{proof}
  Let $V$ be an irreducible $\PP_n$-left module, then either $\ol{V}=V$,
  hence $V\cong_{\PP_n} M^0$, or $\ol{V}=0$. In the latter case, let
  $v\in V\backslash\{0\}$, then $\PP_nv$ is a $\PP_n$-sub-module $\neq
  0$ of $V$.  It follows that the $\PP_n$-left module homomorphism
  $\alpha:\PP_n\to V,\,\varphi\mapsto \varphi v$ is surjective. Standard
  arguments based on Corollary~\ref{rad-lambda} and
  $$
  \PP_n/\Rad\,\PP_n 
  \cong_{\PP_n}
  \bigoplus_{p\podd  n} \Lambda^p/\Rad_{\PP_n}\Lambda^p 
  = 
  \bigoplus_{p\podd  n} M^p
  $$
  complete the proof of the first claim, since $\ker\,\alpha$ is a
  maximal sub-module of $\PP_n$, thus contains $\Rad\,\PP_n$.
  
  Assume now that $q\p n$ contains an even letter.  Let $t\zerl n$
  such that $t\ass q$ and $r\zerlodd n$, then
  $\gamma_r\gamma_t\in\ker\,c\cap\Lambda^q=\Rad_{\DD_n}\Lambda^q$, by
  \kref{c-von-highlieid} and \cite[Corollary~1.6]{blelau96}, hence
  $\gamma_r M^q=0$. This shows $M^q\cong M^0$ (as a $\PP_n$-module).
\end{proof}

A direct consequence of \cite[Proposition~1.1]{blelau02} is

\begin{proposition} \label{cartan-helfer}
  Let $M$ be a $\DD_n$-left module, then, for each $p\podd  n$, the
  multiplicity of $M^p$ in a $\PP_n$-composition series of $M$ is
  $$
  \mlt{M}{M^p}{\PP_n} := \dim\,\gamma_pM,
  $$
  which is also the multiplicity $\mlt{M}{M^p}{\DD_n}$ of $M^p$ in
  a $\DD_n$-composition series of $M$.  The multiplicity of $M^0$ in a
  $\PP_n$-composition series of $M$ is
  $$
  \mlt{M}{M^0}{\PP_n} 
  = 
  \sum_{\substack{p\p n\\[2pt]p\;\mbox{\footnotesize not odd}}}
  \mlt{M}{M^p}{\DD_n} 
  =
  \dim\,M-\sum_{p\podd n}\mlt{M}{M^p}{\PP_n}.
  $$
\end{proposition}

A description of the Cartan invariants 
$$
c_{qp} := \mlt{\Lambda^q}{M^p}{\PP_n} \qquad(q,p\podd  n)
$$
of $\PP_n$ is now immediate from \cite{blelau02}. Recall that any
word $r\in\N^*$ has a unique decomposition $ r = r^{(1)}\wldots
r^{(l)} $ into Lyndon words $r^{(1)},\ldots,r^{(l)}\in\N^*$ such that
$ r^{(1)}\ge_{lex}\cdots\ge_{lex}r^{(l)} $, where $\le_{lex}$ denotes
the lexicographical order on $\N^*$
(\cite[Proposition~7.1]{garsia89}).  If $r$ is a composition of $n$,
then so is the \emph{Lyndon sum composition} corresponding to $r$,
defined by
$$
\LSC(r) := (\summe r^{(1)})\wldots(\summe r^{(l)}).
$$

\begin{corollary} \label{cartan-pn}
  The Cartan matrix of $\PP_n$ is given by
  $$
  c_{qp}
  = 
  \mlt{\Lambda^q}{M^p}{\DD_n} 
  =
  \#\setofall{r\zerl n}{r\ass q,\,\LSC(r)\ass p},
  $$
  for all $p,q\podd  n$.
\end{corollary}

\begin{proof}
  The first equality follows from Proposition~\ref{cartan-helfer}, the
  second from \cite[Corollary~2.1 (and its proof)]{blelau02} (see also
  \cite[Theorem~5.4]{garsia-reutenauer89}), since
  $$
  \mlt{\Lambda^q}{M^p}{\DD_n} 
  =
  \dim\,\gamma_q\Lambda^p
  =
  \dim\,\gamma_q\DD_n\gamma_p
  =
  \dim\,\omega_q\DD_n\omega_p\,,
  $$
  by Propositions~\ref{indec-dn}, \ref{cartan-helfer}, and \kref{eq:pr}.
\end{proof}

Let $A$ be a finite dimensional $K$-algebra, then the \emph{descending
  Loewy series} of an $A$-module $M$ is the chain of sub-modules
$\Rad_A^{(j)}M$, defined by $\Rad_A^{(0)}M:=M$ and
$\Rad_A^{(j)}M:=\Rad_A(\Rad_A^{(j-1)}M)$ for all $j\in\N$.

To conclude this section, enough information on the descending Loewy
series of $\PP_n$ is presented to determine the nil-potency index of
$\Rad\,\PP_n$.

\begin{lemma} \label{loewy-helper}
  Let $M$ be a $\PP_n$-left module and $U$ be a 
  sub-module of $M$ such that
  \begin{listing}
    \item
      $M/U$ is a zero module;
    \item
      $U$ is completely reducible;
    \item
      $\mlt{U}{M^p}{\PP_n}\neq 0\neq \mlt{U}{M^q}{\PP_n}$
      and
      $p\assozerl q$ imply $p=q$, for all $p,q\podd n$;
  \end{listing}
  then $M$ is completely reducible.  
\end{lemma}

\begin{proof}
  Assume first that $U$ has co-dimension $1$ in $M$, and let $m\in M$
  such that $M=U\oplus\erz{m}_K$.  
  Set $\tilde{m}:=m-\sum_{p\podd  n}\frac{1}{p?}\gamma_pm$, then
  $M=U\oplus\erz{\tilde{m}}_K$, since
  $\sum_{p\podd  n}\frac{1}{p?}\gamma_pm\in U$, by (i).

  We will show that $\erz{\tilde{m}}_K$ is a trivial sub-module of $M$, which 
  implies the claim, by (ii).

  Let $r\zerlodd n$ and denote by $q$ the partition obtained by rearranging
  $r$, then one of the following two cases holds.
  
  \textbf{case 1:} $\dim\,\gamma_qM=\mlt{M}{M^q}{\PP_n}=0$, then 
  $\gamma_r\tilde{m}=\frac{1}{q?}\gamma_r\gamma_q\tilde{m}=0$.
  
  \textbf{case 2:} $\dim\,\gamma_qM=\mlt{M}{M^q}{\PP_n}\neq 0$, then
  (iii), \kref{eq:pr} and \kref{eq:lia} imply
  $\gamma_r\gamma_p=\frac{1}{q?}\gamma_r\gamma_q\gamma_p=0$ for all
  $p\neq q$ such that $\mlt{M}{M^p}{\PP_n}\neq 0$, while $\gamma_pm=0$
  for all $p\podd n$ such that $\mlt{M}{M^p}{\PP_n}=0$.  It follows
  that
  $$
  \gamma_r\tilde{m}
  =
  \gamma_rm-\sum_{p\podd  n}\frac{1}{p?}\gamma_r\gamma_pm
  =
  \gamma_rm-\frac{1}{q?}\gamma_r\gamma_qm
  =
  0,
  $$
  by \kref{eq:pr}.

  If the co-dimension $d$ of $U$ in $M$ is larger then $1$, then choose
  a maximal sub-module $\hat{M}$ of $M$ containing $U$.  $\hat{M}$ is
  completely reducible, by induction on $d$, and has co-dimension $1$
  in $M$, by (i). Now apply the part already proven to $\hat{M}$
  instead of $U$ to conclude that $M$ is completely reducible.
\end{proof}

\begin{corollary} \label{rad-inclusion}
  $
  \Rad^{(j)}_{\PP_n}\Lambda^p
  \subseteq 
  \Rad^{(2j-1)}_{\DD_n}\Lambda^p,
  $
  for all $j\in\N$, $p\podd  n$.
\end{corollary}

\begin{proof}
  Let $R^i:=\Rad^{(i)}_{\DD_n}\Lambda^p$ for all $i\in\N_0$, then
  $R^i/R^{i+1}$ is completely reducible as a $\DD_n$-module.
  
  Let $i\in\N_0$ and $r\p n$ such that $M_r$ is isomorphic to an
  irreducible constituent of $R^i/R^{i+1}$, then
  $\length(r)=\length(p)-i$, by \cite[Theorem~2.2]{blelau02}.  In
  particular, for each $j\in\N$, $R^{2j-1}/R^{2j}$ is a zero module
  for $\PP_n$, since $r$ is not odd whenever
  $\length(r)=\length(p)-(2j-1)$.  Furthermore,
  $\length(r)=\length(s)$ for all $r,s\p n$ such that $M^r$ and $M^s$
  occur in $R^{2j}/R^{2j+1}$, hence $r\assozerl s$ implies $r=s$.
  
  Applying Lemma~\ref{loewy-helper} to $M=R^{2j-1}/R^{2j+1}$,
  $U:=R^{2j}/R^{2j+1}$ yields complete reducibility of
  $R^{2j-1}/R^{2j+1}$ as a $\PP_n$-module.  Now
  $\Rad_{\PP_n}\Lambda^p=\Rad_{\DD_n}\Lambda^p$ allows to complete the
  proof by a simple induction.
\end{proof}

We conjecture that the descending Loewy series of $\Lambda^p$ (as a
$\PP_n$-module) is equal to the sub-chain of the descending Loewy
series of $\Lambda^p$ (as a $\DD_n$-module) consisting of $\Lambda^p$
and the odd steps $\Rad^{(2j-1)}_{\DD_n}\Lambda^p$, $j\in\N$.  A
proof, however, can only be given in the case where $p$ is a hook
partition.

\begin{lemma} \label{rad-haken}
  Let $k\in\haken{n-1}$ be odd and $p=k.1^{.(n-k)}\podd  n$,
  then
  $$
  \Rad^{(j)}_{\PP_n}\Lambda^p
  =
  \Rad^{(2j-1)}_{\DD_n}\Lambda^p,
  $$
  for all $j\in\N$.
\end{lemma}

\begin{proof}
  For $j=1$, this has been proven in Corollary~\ref{rad-lambda}.  Let
  $j>1$.  
  
  If $k\ge n-2$, then $\Rad^{(2j-1)}_{\DD_n}\Lambda^p=0$, by
  \cite[Corollary~2.6]{blelau96}, and the claim follows from
  Corollary~\ref{rad-inclusion}.
  
  Let $k<n-2$ and $s\p n$ such that $\length(s)=\length(p)-1$ and
  $p\assozerl s$, then $s=(k+1).1^{.(n-k-1)}$ or $s=k.2.1^{.(n-k-2)}$.
  In the latter case, $\gamma_s\gamma_p=0$, by
  \cite[Proposition~2.1]{blelau96}, since $\gamma_s\gamma_p\in
  \omega_s\DD_n\omega_p\gamma_p$, by \kref{eq:pr}.  Applying
  \cite[Corollary~2.4]{blelau96} and Proposition~\ref{indec-dn} thus
  yields
  $$
  \Rad^{(2j-1)}_{\DD_n}\Lambda^p 
  =
  (\Rad^{(2j-2)}_{\DD_n}\Lambda^r)\omega_r\gamma_p\,,
  $$
  where $r=(k+1).1^{.(n-k-1)}$.  The same argument, applied to $r$
  instead of $p$, gives
  $$
  \Rad^{(2j-2)}_{\DD_n}\Lambda^r 
  =
  (\Rad^{(2j-3)}_{\DD_n}\Lambda^s)\omega_s\gamma_r\,,
  $$
  where $s=(k+2).1^{.(n-k-2)}$.  Furthermore,
  $
  \omega_s\omega_r\gamma_p
  \in
  \ker\,c\cap\Lambda^p
  =
  \Rad_{\PP_n}\Lambda^p
  $, by Theorem~\ref{sol-rad}, \kref{c-von-highlieid} and
  Corollary~\ref{rad-lambda}.  Putting together all pieces,
  \begin{eqnarray*}
    \Rad^{(2j-1)}_{\DD_n}\Lambda^p
    & = &
    (\Rad^{(2(j-1)-1)}_{\DD_n}\Lambda^s)
    \omega_s\gamma_r\omega_r\gamma_q \\[1mm]
    & \subseteq &
    (\Rad^{(j-1)}_{\PP_n}\Lambda^s)\omega_s\omega_r\gamma_q \\[1mm]
    & \subseteq &
    (\Rad^{(j-1)}_{\PP_n}\Lambda^s)\Rad_{\PP_n}\Lambda^q \\[1mm]
    & \subseteq &
    (\Rad^{(j-1)}_{\PP_n}\PP_n)\Rad_{\PP_n}\Lambda^q \\[1mm]
    & \subseteq &
    \Rad^{(j)}_{\PP_n}\Lambda^q,    
  \end{eqnarray*}
  by \kref{eq:pr} and induction.
\end{proof}

\begin{theorem} \label{nilindex}
  The nil-potency index of $\Rad\,\PP_n$ is $\frac{n-1}{2}$, if $n$ is
  odd, and $\frac{n}{2}$, if $n$ is even.
\end{theorem}

\begin{proof}
  If $n=2$, then $\PP_n=\Lambda^{1.1}$ and $\Rad\,\PP_n=0$.  Let $n\ge
  3$. As a consequence of \cite[Corollary~2.6]{blelau96},
  $$
  \Rad_{\DD_n}^{(n-1)}\Lambda^p=\Rad_{\DD_n}^{(n-2)}\Lambda^p=0
  \mbox{ for all $p\podd  n$},
  $$
  since the only partition of $n$ of length $n-1$ is $q=2.1^{.(n-2)}$,
  and $q$ is not odd.
  Corollary~\ref{rad-inclusion} therefore implies that the
  nil-potency index of 
  $$
  \Rad\,\PP_n 
  = 
  \bigoplus_{p\podd n}\Rad_{\PP_n}\Lambda^p, 
  $$
  is less or equal $n/2$. We mention that this part of the claim
  can also be deduced directly from \kref{eq:pr}, \kref{eq:lia} and
  the description of $\Rad\,\PP_n$ given in Corollary~\ref{rad-lambda}.

  Conversely, consider $p=3.1^{.(n-3)}\podd n$, then
  $\Rad_{\DD_n}^{(n-3)}\Lambda^p\neq 0$, by
  \cite[Corollary~2.6]{blelau96} again.  Applying
  Lemma~\ref{rad-haken} yields $\Rad_{\PP_n}^{(n/2-1)}\Lambda^p\neq
  0$, if $n$ is even, and $\Rad_{\PP_n}^{((n-3)/2)}\Lambda^p\neq 0$,
  if $n$ is odd, which completes the proof.
\end{proof}

\newcommand{\appears}[1]{to appear in \emph{#1}}

\end{document}